\documentclass[leqno,12pt]{amsart}

\makeatletter
\renewcommand\paragraph{\@startsection{paragraph}{4}{\z@}%
            {-2.5ex\@plus -1ex \@minus -.25ex}%
            {1.25ex \@plus .25ex}%
            {\normalfont\normalsize\bfseries}}
\makeatother
\makeatletter
\renewcommand{\@seccntformat}[1]{%
  \ifcsname format#1\endcsname
    \csname format#1\endcsname
  \else
    \csname the#1\endcsname\quad
  \fi
}
\makeatother

\usepackage{amsmath}\usepackage{MnSymbol}\usepackage{phaistos}
\usepackage{amsfonts}\usepackage{verbatim}\usepackage{bbold}
\usepackage{graphicx}\usepackage{fontawesome}
\usepackage{marvosym}\usepackage{ifsym}
\usepackage{dictsym}\usepackage{bbm}
\usepackage{wasysym}\usepackage{dutchcal}
\usepackage{amstext}
\usepackage{amsbsy}
\DeclareFontFamily{U}{mathb}{}
\DeclareFontShape{U}{mathb}{m}{n}{
     <-5.5> mathb5
  <5.5-6.5> mathb6
  <6.5-7.5> mathb7
  <7.5-8.5> mathb8
  <8.5-9.5> mathb9
  <9.5-11>  mathb10
  <11->     mathb12
}{}
\DeclareSymbolFont{mathb}{U}{mathb}{m}{n}
\DeclareMathSymbol{\lefttorightarrow}{3}{mathb}{"FC}
\DeclareMathSymbol{\righttoleftarrow}{3}{mathb}{"FD}

\makeatletter
\DeclareRobustCommand{\looparrow}[1]{%
  \mathrel{\mathpalette\looparrow@{#1}}%
}
\newcommand{\looparrow@}[2]{\reflectbox{$\m@th#1#2$}}
\makeatother

\usepackage{amsopn}\usepackage[
    colorlinks=false,
    linkcolor=blue,
    bookmarks=true,
    bookmarksopen=true,
    bookmarksnumbered=true, pdfpagelabels
]{hyperref}
\usepackage{amsthm}\usepackage{color}
\usepackage{pigpen}
\theoremstyle{definition}

\theoremstyle{remark}

\def\proclaim#1{\vskip0.5em\noindent{\bf #1}\it}
\def\endproclaim{\vskip0.5em\par\noindent\rm}
\def\proclaim#1{\vskip0.5em\noindent{\bf #1}\it}
\def\endproclaim{\vskip0.5em\par\noindent\rm}
\def\demo#1{\vskip0.5em\noindent{\bf #1\ }}

\def\text#1{\mbox{#1}}
\def\flushpar{\par\noindent}

\def\tag#1{\eqno{(#1)}}

\def\e{\varepsilon}
\def\a{\alpha}
\def\b{\beta}
\def\G{\Gamma}\def\g{\gamma}
\def\d{\delta}
\def\th{\theta}
\def\D{\Delta}
\def\s{\sigma}
\def\l{\lambda}
\def\x{\times}

\def\o{\overline}
\def\f{\flushpar}
\def\u{\underline}
\def\v{\varphi}

\def\om{\omega}

\def\B{\mathcal B}

\def\T{\widehat T}
\def\({\biggl(}
\def\){\biggr)}
\def\<{\langle}
\def\>{\rangle}
\def\bdy{\partial}\def\({\biggl(}
\def\){\biggr)}
\def\[{\biggl[}
\def\]{\biggr]}
\def\bul{\smallskip\f$\bullet\ \ \
$}\def\sbul{\f$\bullet\ \ \ $}\def\Par{\smallskip\f\P}
\def\pf{\smallskip\f{\tt Proof}\ \ \ \ }\def\sms{\smallskip\f}
 \def\lfl{\lfloor}\def\rfl{\rfloor}\def\st{\text{such that}}
\def\ttau{\widehat{\tau}}

\let\oldtocsection=\tocsection

\let\oldtocsubsection=\tocsubsection

\let\oldtocsubsubsection=\tocsubsubsection

\renewcommand{\tocsection}[2]{\hspace{0em}\oldtocsection{#1}{#2}}
\renewcommand{\tocsubsection}[2]{\hspace{1em}\oldtocsubsection{#1}{#2}}
\renewcommand{\tocsubsubsection}[2]{\hspace{2em}\oldtocsubsubsection{#1}{#2}}
\begin{document}
\title{Dynamics of inner functions   revisited}

\author{Jon Aaronson and Mahendra Nadkarni}
\address[Aaronson]{School of Math. Sciences, Tel Aviv University,
69978 Tel Aviv, Israel.}
\email{aaro@tau.ac.il}
\address[Nadkarni]{Department of Mathematics, University of Mumbai, Kalina, Mumbai 400098, India}
\email{mgnadkarni@gmail.com}\begin{abstract}

 We study the circle restrictions of inner functions of the unit disc showing that the local invertibility of a restriction is independent of its singularity set and proving a local characterization of analytic conditional expectations.

 We establish   central limit properties for some stochastic processes driven by probability preserving restrictions via spectral analysis of their perturbed transfer operators.
 \end{abstract}
\thanks{\copyright 2022-25.}
\maketitle\markboth{\copyright J. Aaronson and M.  Nadkarni}{Inner functions}
\tableofcontents
\setcounter{tocdepth}{0}
\section{ Introduction}

\

\subsection{Inner functions and their restrictions}
\

An {\it inner function} of the {\it  unit disc} $\mathbb{D}:=\{z\in\Bbb C:\ |z|<1\}$ is an analytic endomorphism $\phi:\mathbb{D}\righttoleftarrow$ so that
for Lebesgue almost every $\xi\in\bdy\mathbb{D}$,
$$\phi (r\xi)\xrightarrow[r\to 1-]{×} \phi(\xi)\in\bdy\mathbb{D}$$
The {\it
restriction} of $\phi$ to $\bdy\mathbb{D}$ is defined $\l$-a.s. (where $\l$ is Lebesgue measure on $\bdy\mathbb{D}$) and is
  a {\it nonsingular transformation} of
$(\bdy\mathbb{D},\l)$ in the sense  that  $\l$ and $\l\circ \phi^{-1}$ have the same null sets.

\

This follows from Nordgren's theorem (in \S\ref{StructRestr}) which also shows the connection between the ergodic theory of an inner function restriction on $(\bdy\mathbb{D},\l)$ and the dynamics of the action of the inner function on $\mathbb{D}$.

For this, and more discussion of the ergodic theory of restrictions, see e.g. \cite[Ch.6]{A1} and references therein.
\

\subsection{Overview of the paper}
This paper deals with the structure and properties of inner functions, the spectral theory of their transfer operators and the central limit theory of stochastic processes driven by their restrictions.

The rest of this subsection
is devoted to a description of the main results.

\subsubsection*{Local invertibility vs. singularity set}
\

 In \S\ref{loc-invert} we consider  {\tt\small local invertibility}   of general nonsingular maps
 showing {\it i.a.} that  a nonsingular transformation of a standard, nonatomic probability space is  locally invertible    iff it is
\f {\tt\small forward nonsingular}  (Theorem 2.2,\ p.\pageref{struct}).
 \

 This  enables
 an elementary   proof of a multiplicity result of   Aleksandrov
 (Proposition 2.7,\ on p.\pageref{AH}):  an inner function has a locally invertible restriction iff it admits angular derivatives a.e. on $\bdy\mathbb{D}$.
\

The {\tt singularity set} of an inner function
\footnote{the (closed) set of points of $\bdy\mathbb{D}$ where it is not  analytic} (see \S\ref{sing})
was studied in
in \cite{Seidel}.
Evidently, if an inner function has a Lebesgue-null singularity set, then it's
restriction, being analytic a.e. on $\Bbb T$, is locally invertible. However (by Proposition 2.8 on p.\pageref{prop2.7}),    any closed set of $\bdy\mathbb D$ appears as the singularity set of an inner function  whose restriction is locally invertible.

\subsubsection*{Spectrum  of the transfer operator and central limits}
\

 In \S\ref{QcTO} $\&$ \S\ref{pert}, we restrict attention to  non-M\"obius, inner functions $\phi:\mathbb{D}\righttoleftarrow$  preserving an absolutely continuous probability; showing in \S\ref{QcTO} that
their transfer operators have spectral gaps\footnote{  i.e.  are Doeblin-Fortet operators (as on p.\pageref{DF(i)})} on  weighted Hilbert spaces (see \S\ref{WHS}).

In \S\ref{pert} we consider  central limit properties of stochastic processes
$(\psi\circ\tau^n:\ n\ge 1)$  ($\psi:\Bbb T\to\Bbb R^d$) driven by such restrictions
  $\tau=\tau(\phi)$.    In particular, if  $\phi:\mathbb{D}\righttoleftarrow$ has a singularity on $\bdy\mathbb{D}$ and $\psi:\mathbb{T}\to\mathbb{R}$ is   real analytic and not constant, then $(\psi\circ\tau^n:\ n\ge 1)$ satisfies the { conditional central limit theorem} (as in \eqref{con-CLT}
  on p.\pageref{con-CLT}).   These results rely on smooth perturbations of the quasicompact transfer operators under consideration.
 \

 The {central limit theorem} for measure preserving inner function restrictions
  (as in \eqref{CLT}   on p.\pageref{CLT}) is established in \cite{OlegMar} for e.g. non-constant H\"older continuous functions.
 See also \cite{N-S} for  different forms of central limit theorem for measure preserving inner function restrictions.

 \subsubsection*{Analytic conditional expectations}
 \

 Aleksandrov proved in \cite{Aleksandrov86} that conditional expectation with respect to a sub-$\s$-algebra $\mathcal{C}\subset\B(\Bbb T)$ is {\it analytic} in the sense that the conditional expectation projection commutes with the {\tt\small Riesz projection} (orthogonal projection $L^2\to H_0^2:=\{f\in H^2\ \hat{f}(0)=0\}$ )  if and only if $\mathcal{C}=\tau^{-1}\B(\Bbb T)$
where $\tau=\tau(\phi)$ with $\phi$ inner,\ $\phi(0)=0$.
\

 In \S\ref{loc-Aleks} we prove a local version of this  (Theorem 5.2 on p. \pageref{thm5.2}).
 \

 \

  \section{ Structure and basic properties of inner functions}
 
   \subsection{ Nonsingular maps and transformations}\label{NST}
   \

  A {\it nonsingular  map} $\pi:(X,m)\to (Y,\mu)$  of the non-atomic, Polish probability spaces $(X,m)\ \&\ (Y,\mu)$  is a
 measurable map $\pi:X\to Y$  satisfying $m\circ \pi^{-1}\sim \mu$. It is {\it probability preserving} if $m\circ T^{-1}= \mu$.

 A {\it nonsingular}  [probability preserving] {\it transformation}   $(X,m,T)$ is a
 nonsingular [{\emph probability preserving}] self map $T:(X,m)\to (X,m)$.
 \
 

        \
        
        The {\it transfer operator} (aka ``adjoint'' or ``predual'') of the nonsingular map $\pi:(X,m)\to (Y,\mu)$  is $\widehat{\pi}:L^1(m)\to L^1(\mu)$ defined  by
        $$\int_A\widehat{\pi} fd\mu\ =\ \int_{\pi^{-1}A} fdm\ \ \ A\in\B(X).$$
     
           \
           
  The following is a standard application of the disintegration theorem (\cite[Ch. III]{De-Meyer}, also \cite[Ch.1]{A1}):
   \

   \proclaim{Proposition 2.1\ \ (Preimage measures)}
   \
   
 Let  $\pi:(X,m)\to (Y,\mu)$ be  a  nonsingular map,
 then
 \

           $\exists\ Y_0\in\B(X),\ \mu(Y_0)=1$ and
           $\nu=\nu^{(\pi)}:Y_0\to\mathcal{M}(X)$ so that
           \begin{align*}\tag*{\dsjuridical}\label{dsjuridical}&\widehat{\pi}  1_A(y)=\nu_y(A)\ \text{for}\ y\in Y_0,\ A\in\B(X)\ \&\ \nu_y(X)=\tfrac{dm\circ \pi^{-1}}{d\mu}(x).\\ &
           \text{Moreover,} \  \nu_y(X\setminus \pi^{-1}\{y\})=0\ \&\ \nu_x\perp m\ \text{for $\mu$-a.e.}\ y\in Y_0.
           \end{align*}

  \endproclaim
      The function $x\mapsto\nu_x$ is known as the {\it transition kernel} of $\widehat{\pi} $ and the $\nu_x$      are known as {\it preimage -} or {\it fiber measures}.

       \subsection{ Local invertibility and forward nonsingularity}\label{loc-invert}
        \
        
    We'll say that the nonsingular map $\pi:(X,m)\to (Y,\mu)$ is
      \sbul {\it almost countable to one} if $\exists\ Y_0\in\B(Y),\ \mu(Y\setminus Y_0)=0$
        so that $\pi^{-1}\{y\}\cap Y_0$ is at most countable $\forall\ y\in Y_0$;
        \sbul {\it locally invertible} if $\exists$ an at most countable partition $\a\subset\B(Y)$
      so that    $\pi:a\to \pi a$ is invertible, nonsingular $\forall\ a\in\a$;
              \sbul {\it forward nonsingular} if $\exists\ X_0\in\B(X),\ m(X\setminus X_0)=0$ so that
$A\in\B(X_0),\ m(A)=0\ \implies\ \mu(\pi A)=0.$
    
    As shown in \cite{Rokhlin} (see also \cite{Nad-inner} and \cite[Ch. 1]{A1}), for the nonsingular map $\pi:(X,m)\to (Y,\mu)$ the conditions almost countable to one and  locally invertibility are both equivalent to the pure atomicity of its transition kernel (i.e. almost every preimage measure is purely atomic).

\proclaim{Theorem 2.2}\label{struct}
\
The nonsingular nonsingular map $\pi:(X,m)\to (Y,\mu)$ is forward nonsingular iff it is locally invertible.

\endproclaim
\demo{Proof}

    It is standard that local invertibility implies forward nonsingularity.
    
    We``ll show that forward nonsingularity $\implies$ pure atomicity of the transition kernel.
    \
    
    Let $y\mapsto \nu_y$ ($Y_0\to\mathcal{M}(X)$) be the transtition kernel of $\widehat{\pi}$.
    \

    We claim first that the function $x\mapsto\nu_{\pi x}(\{x\})$ is measurable.
    \
    
    To see this let
    $d$ be a Polish metric on $X$ and let $\a_n$ be a sequence of partitions of $X$ so that
    $\a_{n+1}\succ\a_n$,\ 
    $\sup_{a\in\a_n}\,\text{\tt diam}\, a\xrightarrow[n\to\infty]{}0$.
    \

    Thus $\a_n(y)\downarrow\ \{y\}\ \forall\ y\in X$ (where
    $y\in\a_n(y)\in\a_n$) and the function $y\mapsto\nu_{\pi y}(\{y\})=\lim_{n\to\infty}\nu_{\pi y}(\a_n(y))$  is measurable.
    
    It follows that $W:=\{x\in X:\ \nu_{\pi x}(\{x\})=0\}\in\B(X)$.
    \

    If the transition kernel is not purely atomic, then $m(W)>0$.
    \

Next, $\pi W$ is analytic, whence  universally measurable and
since $\pi ^{-1}\pi W\supseteq W$,\ we have $m(\pi^{-1}\pi W)\ge m(W)>0$ whence (by nonsingularity of $\pi$), $\mu(\pi W)>0$.
\

For $V\in\B(W)$,
\begin{align*}
 m(V) &=m(V\cap W)=m(V\cap \pi^{-1}\pi W)\ \ \because\ \pi^{-1}\pi W\supset W\\ &=
 \int_{\pi W}\widehat{\pi} 1_Vd\mu=\int_{\pi W}\nu_x(V)d\mu(x).
\end{align*}


By the analytic section theorem (\cite{Lusin}, see also \cite{Jankoff}, \cite{von-N-section}
\cite[Thm 5.5.2]{Srivastava}, \cite[\S 8.5]{Coh}); and Lusin's continuity theorem,
\

$\exists\ A\subset\ \pi W$  compact,\ $\mu(A)>0$ , and $\zeta:A\to B:=\zeta A\subset W$ continuous so that $\pi \circ \zeta=\text{\tt\small Id}$.
   \
   
   It follows that 
   $$B\cap \pi^{-1}\{x\}=\begin{cases}&\emptyset\ \ \ \ \ x\notin A;
   \\ &\{\zeta(x)\}\ \ \ \ \ x\in A
                       \end{cases}$$
                     whence
$$\nu_x(B)=1_A(x)\nu_x(\{\zeta(x)\})=0\ \because\ \pi(x)\in W\ \&\ x=\pi(\zeta(x)).$$
Thus
   \begin{align*}
    m(B)=\int_A\nu_x(\{\zeta(x)\})d\mu(x)=0
   \end{align*}
whereas $\pi B=\pi \zeta A=A,\ \mu(\pi B)=\mu(A)>0$ and forward nonsingularity fails.\ \CheckedBox

  \subsection{Structure of inner function restrictions}\label{StructRestr}
\

It will be convenient to identify $\bdy\mathbb{D}$ with
$\mathbb{T}:=\mathbb{R}/\mathbb{Z}\cong [0,1)$\
via\ \ $$ x\in [0,1)\leftrightarrow\ \chi(t):=e^{2\pi ix}\in\bdy\mathbb{D}.$$
Indeed, $\chi:(\Bbb T,m)\to (\bdy\mathbb{D},\l)$ is an isomorphism of measure spaces where $m$ is Lebesgue measure on $[0,1)$.
\

Let $\phi:\mathbb{D}\righttoleftarrow$ be inner with
restriction $\phi:\bdy\mathbb{D}\righttoleftarrow$. Define $\tau=\tau(\phi ):\Bbb T\righttoleftarrow$  $m$-a.e. by $\tau(t):= \chi^{-1}(\v(\chi(t)))$.

\

The following shows that $(\Bbb T,m,\tau)$ is a nonsingular transformation, whence so is the isomorphic $(\bdy\mathbb{D},\l,\v)$.

\proclaim{Nordgren's Theorem \ \cite{Nor}}
\

For $\phi:\mathbb{D}\righttoleftarrow$ inner,
\begin{align*}\tag*{\dsheraldical}\label{dsheraldical}
 \pi_z\circ \tau^{-1}=\pi_{\phi (z)},\ \text{\rm\footnotesize or,  equivalently:}\ \ttau p_z=p_{\phi(z)}.
\end{align*}
where $\pi_z$\ {\rm[$p_z$]} is the {\tt\small Poisson measure [density] at $z$} defined by
\begin{align*}\tag*{\PHtunny}\label{PHtunny}
 d\pi_z(x):=p_z(x)dx,\ p_z(x):=\text{\tt Re}\tfrac{\chi(x)+z}{\chi(x)-z}.
\end{align*}
 \endproclaim
 \f{\bf Remark:}\ \ Equation \eqref{dsheraldical} \ is aka {\tt\small Boole's formula} as a version  for rational inner functions of the  the upper half plane
 appears in \cite[p. 787]{Boo} (see also \cite[\S8]{Glaisher}).   For a converse to Nordgren's theorem,  see  \cite{Let}.
 \

 \subsection{ Invariant probabilities for inner function restrictions}
\

 \proclaim{Denjoy-Wolff Theorem\ \ \cite{Denjoy,Wolff-I}}
\par Suppose that $\phi:\mathbb{D}\righttoleftarrow$ is analytic, not constant and not M\"obius, then
there is a (unique) point $\mathfrak{d}=\mathfrak{d}(\phi)\in\o{\mathbb{D}}$ such that
$${1-|\phi(z)|^2\over |1-\o{\mathfrak{d}}\phi(z)|^2}
\ge {1-|z|^2\over |1-\o{\mathfrak{d}}z|^2}\ \forall\ z\in \mathbb{D},$$
and
$$\phi^n(z)\to\mathfrak{d}\text{ as }n\to\infty\ \ \forall\ z\in \mathbb{D}.$$
\endproclaim
The point $\mathfrak{d}(\phi)$ is called the {\it Denjoy-Wolff} point of $\phi$.
                                                                
 \proclaim{Corollary 2.3\ \  ({\rm see e.g. \cite{A-inner78, Neu, Man-Do}})}

 The restriction $\tau$ of an inner function $\phi$ has an absolutely continuous, invariant probability iff the
Denjoy-Wolff point of $\phi$ belongs to $\mathbb{D}$ and in this case $(\Bbb T,\pi_{\mathfrak{d}(\phi)},T)$ is either conjugate to a circle rotation or is an exact probability preserving transformation.\endproclaim

 \subsubsection{Clark measures}\label{Clark}
 \

As in \cite{Saksman}, the {\it Clark measure} of an analytic endomorphism
 $\phi:\Bbb D\righttoleftarrow$ at $\xi\in\bdy\Bbb D$ is the representing measure $\mu_\xi=\mu_\xi^{(\phi)}\in\mathcal{M}(\Bbb T)$ of the positive harmonic function
 $z\mapsto \text{\tt Re}\,\tfrac{\xi+\phi(z)}{\xi-\phi(z)}$   satisfying
 \begin{align*}\tag*{\Bat}\label{Bat}
  \text{\tt Re}\,\tfrac{\xi+\phi(z)}{\xi-\phi(z)}=\int_{\Bbb T}p_zd\mu_\xi.
 \end{align*}
It follows that $\xi\mapsto\mu_\xi$ is weak $*$ continuous ($\bdy\Bbb D\to\mathcal{M}(\Bbb T)$),
\

Now let $\phi:\Bbb D\righttoleftarrow$ be inner,
then $\text{\tt Re}\,\tfrac{\xi+\phi(r\chi)}{\xi-\phi(r\chi)}\xrightarrow[r\to 1-]{}0$ a.s. whereas by Fatou's theorem
$\int_{\Bbb T}p_{r\chi}d\mu_\xi\xrightarrow[r\to 1-]{}\tfrac{d\mu_\xi}{dm}$ a.s and $\mu_\xi\perp m\ \forall\ \xi\in\bdy\Bbb D$.
 \

 Moreover, since $p_z(t)=\text{\tt Re}\,\tfrac{\chi(t)+z}{\chi(t)-z}$, it follows from \eqref{dsheraldical} (on p. \pageref{dsheraldical}) that $\nu_t^{(\tau)}=\mu_{\chi(t)}$ where
 $\nu^{(\tau)}:\Bbb T\to\mathcal{M}(\Bbb T)$ are the preimage measures of the restriction $\tau=\tau(\phi)$. See also \cite[\S2]{Saksman}.

 \subsection{Factorization of  inner functions}

 \subsubsection{ Blaschke products}
\

Suppose that $Z\subset\mathbb{D}$ is countable and that $\mathfrak{m}:Z\to\Bbb N$ is so that
 $\sum_{a
 \in Z} \mathfrak{m}(a)(1-|a|)<\infty$. 
 \
 
 The  {\it Blaschke product} with zero set $Z$ and {\it multiplicity function} $\mathfrak{m}$ is 
 $B=B_{Z,\mathfrak{m}}:\mathbb{D}\to\Bbb C$ defined by
 \begin{align*}
 B(z):=\prod_{\a\in Z}b_\a(z)^{\mathfrak{m}(\a)}&\ \text{with} \ b_\a(z):=c_\a\tfrac{z-\a}{1-\o{\a}z}\\ & \text{where}\ c_\a=\tfrac{-\o{\a}}{\a}\ (\a\ne 0)\ \&\ c_0=1.  
 \end{align*}
 Te product converges locally uniformly on $\mathbb{D}$ because
 $$|1-b_\a(z)^{\mathfrak{m}(\a)}|\le\tfrac{1+|z|}{1-|z|}\cdot \mathfrak{m}(\a)(1-|\a|).$$
 \

 It can be shown that $B$ is inner, that $\{a\in\mathbb{D}:\ B(a)=0\}=Z$ and that
 for $a\in Z,\  \ \tfrac{B(z)}{b_a(z)^{\mathfrak{m}(a)}}=:H(z)$ is bounded, analytic with $H(a)\ne 0$.
 \subsubsection{ Singular inner functions}
 \

 An inner function  $S:\mathbb{D}\righttoleftarrow$ is {\it singular} (i.e. without zeros) iff
    $$\log S(z)=-\int_\Bbb T\tfrac{\chi(t)+z}{\chi(t)-z}d\s(t)$$ where $\s\in\frak M(\Bbb T),\ \s\perp m$.
    \
    
    In this situation, we'll denote $S=S_\s$ ($\&$/or $\s=\s_S$).

\
\proclaim{ Factorization Theorem\ \ \cite{Smirnov} ({\rm\small see also \cite[theorem 17.15]{Rudin})}}
\

 Let $\phi:\mathbb{D}\righttoleftarrow$ be inner function,  let
$\mathcal{Z}(\phi):=\{z\in\mathbb{D}:\ \phi(z)=0\}$ and, for\ $a\in\mathcal{Z}(\phi)$, let
$$\mathfrak{m}(a):=\max\,\{k\ge 0:\tfrac{\phi}{b_a^k}\ \text{is bounded on}\ \mathbb{D}\},$$

then   $\sum_{\a
 \in \mathcal{Z}(\phi)} \mathfrak{m}(a)(1-|\a|)<\infty$ and 
$\phi=\l B_{\mathcal{Z}(\phi),\mathfrak{m}}\cdot S$ where $\l\in\Bbb S^1$ and $S$ is a singular inner function.\endproclaim
\

\subsection{Regular points and singularities}\label{sing}
\

A {\it regularity point} of the inner function $\phi:\mathbb{D}\righttoleftarrow$ is a point
$z\in\bdy\mathbb{D}$ so that $\exists\ U\subset\Bbb C$  open with $ z\in U$, and an analytic function $F:U\to\Bbb C$  so that $F|_{U\cap\mathbb{D}}\equiv \phi|_{U\cap\mathbb{D}}$.
\

A non-regular point  in $\bdy\mathbb{D}$ is called a {\it singularity}.
\

The {\it regularity set} of  $\phi$ is $\mathfrak{r}_\phi:=\{\text{\tt\small regularity points of}\ \phi\}$ and the {\it singularity set} of $\phi$ is 
$\frak s_\phi:=\bdy\mathbb{D}\setminus\frak r_\phi.$
\

\bul If $\nu_x$ is any Clark measure for the inner function $\phi$ (as in \S\ref{Clark})
then  $\chi^{-1}\mathfrak{s}_\phi=(\text{\tt spt}\,\nu_x)'$ and
\sbul if $S=S_\s$ is a singular inner function then
$\chi^{-1}\mathfrak{s}_S=\text{\tt spt}\,\s$.

\subsubsection{ Derivative of a restriction at a regular point}
\

Let $\phi:\mathbb{D}\righttoleftarrow$ with restriction $\tau=\tau(\phi)$
, then $\tau$ is differentiable at each $\th\in\chi^{-1}\frak r_\phi$ with
\begin{align*}\tag*{\Industry}\label{Industry}
\begin{split}
 \tau'(\th)=\bdy \phi(\chi(\th))=|\phi'(\chi(\th))|\ \text{with}\ \bdy \phi(z):=\tfrac{z\phi'(z)}{\phi(z)}=z(\log \phi)'(z).
 \end{split}
\end{align*}
If
$\phi=\l B_{\mathcal{Z}(\phi),\mathfrak{m}}S_\s$, then
$$\mathfrak{s}_\phi=\mathcal{Z}(\phi)'\cup\text{\tt spt}\,\s$$
and, as in \cite{Martin-rest}
\begin{align*}\tag*{\faMobile}\label{faMobile}
\begin{split}
\tau'(\th)&= \sum_{\a\in\mathcal{Z}(\phi)}\mathfrak{m}(\a)p_\a(\th)+\tfrac12\int_\Bbb T\tfrac{d\s(t)}{\sin^2(\pi(\th-t))}\\ &\ge \sum_{\a\in\mathcal{Z}(\phi)}\tfrac{1-|\a|}{1+|\a|}+\tfrac{\s(\Bbb T)}2=:\eta\ge 0
\ \ \forall\ \th\in\Bbb T\ \st\ \chi(\th)\in\mathfrak{r}_\phi
\end{split}
\end{align*}
where $p_\a(\th):=\text{\tt Re}\tfrac{\a+\chi(t)}{\a-\chi(t)}$ and $\sum_{\a\in\emptyset}:=0$. Since either $\mathcal{Z}(\phi)\ne\emptyset$ or $\s\ne 0$ (or both), we have $\eta>0$.

\subsubsection{ Arc maps}
\

An {\it arc map} is a 
triple
$(\Bbb T,T,\a)$ where $\a$ a finite or countable
partition $\mod m$ of $\Bbb T$ into {\it open arcs}, (open, connected subsets $A\subsetneq\Bbb T$)\ \ \footnote{of form $\chi^{-1}A=(a,b)$ with $0\le a< b\le 1$ or
$[0,a)\cup (b,1]$  with $0<a\le b<1$};
and $T:\Bbb T\to \Bbb T$ is a map such that
\sbul for each $A\in\a$.
 $T:A\to T(A)$ is a bi-absolutely continuous homeomorphism
 and
 \sbul  $\s(\bigcup_{n\ge 0}T^{-n}\a)=\B(\Bbb T)$.
 \

 It is called {\it piecewise} $C^k$ ($k\ge 1$),
[{\it analytic}] if each  $T:A\to TA$ is a $C^k$-diffeomorphism,
[bi-analytic]  (respectively).
 \

\proclaim{Proposition 2.5\ \ (Arc map restrictions)}\label{prop2.9}
\

Let $\phi:\mathbb{D}\righttoleftarrow$ be inner with $m(\mathfrak{s}_\phi)=0$
and the
Denjoy-Wolff point of $\phi$ belongs to $\mathbb{D}$, then
\sms {\rm (i)}\ \ \ $\tau:\Bbb T\righttoleftarrow$ defines a piecewise
 analytic, arc map $(\Bbb T,\tau,\a)$ such that
 \begin{align*}\tag{U}\label{Exp}
\begin{split}
& \tau^{\prime}(x)\ge\eta>0\ \forall\ x\in \chi^{-1}\mathfrak{r}_\phi\ \&\\ &
\exists\ d\ge 1\ \text{s.t.}\ |\tau^{d\prime}(x)|\ge \b>1\ \ \forall\ x\in \bigcap_{j=0}^{d-1}\tau^{-j}\chi^{-1}\mathfrak{r}_\phi.
 \end{split}
\end{align*}
\sms{\rm (ii)} \ \ If, in addition, $\#\mathcal{Z}(\phi),\ \#\mathfrak{s}_\phi<\infty$, the partition
$\a$ may be chosen to be {\rm surjective}:\begin{align*}\tag{{\tt onto}}\label{onto}
\tau(A)=\Bbb T\ \mod m\ \forall\ A\in\a.
\end{align*}
Moreover
\begin{align*}\tag{A}\label{Adler}
\sup_{x\in \chi^{-1}\mathfrak{r}_\phi}\D\tau(x)<\infty\ \text{with}\
\D\tau(x):=\tfrac{|\tau^{\prime\prime}(x)|}{\tau'(x)^2}.
\end{align*}\endproclaim
 \

We'll call   piecewise onto arc maps satisfying \eqref{Exp},\
 \eqref{onto} and \eqref{Adler} \ {\it Adler maps}.

Adler interval maps (Adler arc maps with surjective partitions into  intervals) are considered in \cite{Fexp}.

Any Adler arc map  is conjugate by rotation to an Adler  interval map.

\

 \endproclaim
\demo{Proof of Proposition 2.5}%
\demo{Proof of \eqref{Exp}}\ \ There is a M\"obius transformation $\Psi:\mathbb{D}\righttoleftarrow$ so that
$\mathfrak{d}(g)=0$ where $g=\Psi^{-1}\circ\phi\circ\Psi$. If the restrictions of $g\ \&\ \Psi$ are $U\ \&\ \psi$ respectively, then $\tau^n=\psi\circ U^n\circ\psi^{-1}\ \forall \ n\ge 1$.
\

Since $g(0)=0$,  by \eqref{faMobile} on p. \pageref{faMobile}, we have
$$U'(\th)=\bdy g(\chi(\th))=1+\bdy h(\chi(\th))\ge 1+\d_h=:\rho>1.$$
For $n\ge 1$,
\begin{align*}
\tau^{n^\prime}&=\psi'\circ U^n\circ\psi^{-1}\cdot U^{n^\prime}\circ\psi^{-1}\cdot \psi^{-1^\prime}\\ &\ge
\rho^n\min_\Bbb T\psi'\cdot\min_\Bbb T\psi^{-1^\prime}\\ &\ge B>1\ \text{for large enough}\ n\ge 1.\ \CheckedBox\ \ \text{\eqref{Exp}}
 \end{align*}
\
\demo{Proof that $\s(\bigcup_{n\ge 0}T^{-n}\a)=\B(\Bbb T)$}\
\

By \eqref{Exp} that $\a_n:=\bigvee_{k=0}^{n-1}\tau^{-1}\a$ is also
a partition $\mod 0$ of $\Bbb T$ into open arcs satisfying
$$\max\,\{m(a):\ a\in\a_n\}\le(\tfrac1\eta)^d\cdot(\tfrac1\b)^{\frac{n}d},$$
whence $\s(\bigcup_{n\ge 0}\tau^{-n}\a)\overset{m}=\B(\Bbb T)$.\ \CheckedBox
\demo{Proof of \eqref{onto} in case $\mathfrak{s}_\phi\ne\emptyset$} \ We construct $\a$, a $\mod 0$ partition of $\Bbb T$ into open arcs satisfying
\eqref{onto}.
\

Since $\#\mathfrak{s}_\phi<\infty$,  \eqref{faMobile}
on p. \pageref{faMobile} now has the form
 \begin{align*}\tag*{\faRocket}\label{faRocket}
\mathfrak{T}'(\th)=\sum_{a\in\mathcal{Z}(\phi)}\mathfrak{m}(a)p_a(\th)+\tfrac12
\sum_{t\in\mathfrak{s}_\phi}\tfrac{\s(\{t\})}{\sin^2(\pi(\th-t))},
   \end{align*}
Suppose that $J\subset\Bbb T$ is an open arc and $f:J\to\Bbb T$ is continuously differentiable on $J$ with $\min_Jf'>0$, then $f$ has a {\it lifting}:
\sbul $\exists$ an interval $\widetilde{J}\subset\Bbb R$ so that $\mathcal{m}\widetilde{J}=J$ where $\mathcal{m}:\Bbb R\to\Bbb T,\ \mathcal{m}(x)=x\ \mod 1$;  and $\exists\ F:\widetilde{J}\to\Bbb R$ continuously differentiable so that
$$\mathcal{m}(F(x))=f(\mathcal{m}(x))\ \text{for}\ x\in \widetilde{J}.$$
In particular, $F'(x)=f'(\mathcal{m}(x))$.
\

Let $J\subset\chi^{-1}\mathfrak{r}_\phi$ be a maximal open arc
(i.e. $\bdy J\subset\chi^{-1}\mathfrak{s}_\phi$)
and let $\mathfrak{T}:\widetilde{J}\to\Bbb R$ be the lifting of $\tau:J\to\Bbb T$.
\

Write $\widetilde{J}=(a_-,a_+)$, then by \eqref{faRocket},

$$\mathfrak{T}'(\th)\xrightarrow[\th\to \{a_-,a_+\},\ \th\in (a_-,a_+)]{}\infty$$
whence
\begin{align*}\tag*{\faMedkit}\label{faMedkit}
\mathfrak{T}(\th)\xrightarrow[\th\to a_\pm,\ \th\in (a_-,a_+)]{}\pm\infty
\end{align*}

and there is a countable $\mod 0$ partition $\mathcal{p}_{\widetilde{J}}$ of $\widetilde{J}$ into open arcs so that for each $A\in\mathcal{p}_{\widetilde{J}}$,  $\mathfrak{T}A$ is an interval of length $1$.
\

It follows that
$\a_J:=\mathcal{m}\mathcal{p}_{\widetilde{J}}$ is a $mod\ 0$ partition of $J$ into open arcs  so that
$\tau A=\Bbb T\ \mod 0\ \forall\ A\in\a_J$.
\

Since $\#\mathfrak{s}_\phi<\infty$ we have that $\chi^{-1}\mathfrak{r}_\phi$ is a finite union of maximal open arcs as above
and so there is a $\mod 0$ partition $\a$ of $\chi^{-1}\mathfrak{r}_\phi$ into open arcs so that $\tau A=\Bbb T\ \mod 0\ \forall\ A\in\a$.\ \CheckedBox\ \eqref{onto}

 \demo{Proof of  \eqref{Adler}}\  Since $\#\mathcal{Z}(\phi)\ \&\ \mathfrak{s}_\phi$ are both finite,  by \eqref{faMobile} on p. \pageref{faMobile},
 $\tau(\th)=\mathfrak{T}(\th)\ \mod 1$ with
 \begin{align*}\tag*{\faPlane}\label{faPlane}
\mathfrak{T}(\th)=
\mathfrak{b}(\th)- \tfrac1{2\pi}\sum_{t\in \mathfrak{s}_\phi}\s(\{t\})\cot(\pi(\th-t))
   \end{align*}
where $\mathfrak{b}\equiv 0$ when $\mathcal{Z}(\phi)=\emptyset$; and
when $1\le \#\mathcal{Z}(\phi)<\infty$;
$$\mathfrak{b}(\th)=\int_0^\th(
\sum_{a\in\mathcal{Z}(\phi)}\mathfrak{m}(a)p_a(t))dt\ \mod 1$$ defines an analytic endomorphism of $\Bbb T$.


\
Since either $\mathfrak{b}=0=\D(\mathfrak{b})$ or $\mathfrak{b}:\Bbb T\righttoleftarrow$
is analytic, expanding in which case $\D(\mathfrak{b}):\Bbb T\to\Bbb R$ is analytic, we have
$\|\D(\mathfrak{b})\|_\infty<\infty$.
\

In case $\mathfrak{s}_\phi\ne\emptyset$, by \eqref{faRocket}
\begin{align*}\tag*{\dsmedical}\label{dsaeronautical}
\tau'(\th)=\mathfrak{b}'(\th)+\tfrac12
\sum_{t\in\mathfrak{s}_\phi}\tfrac{\s(\{t\})}{\sin^2(\pi(\th-t))}=:\mathfrak{b}'(\th)+s'(\th)
\end{align*}
and
\begin{align*}
 |\tau^{\prime\prime}(\th)|&\le|\mathfrak{b}^{\prime\prime}(\th)|+|s^{\prime\prime}(\th)|\\ &\le|\mathfrak{b}^{\prime\prime}(\th)|+\pi
\sum_{t\in\mathfrak{s}_\phi}\s(\{t\})|\tfrac{\cos(\pi(\th-t))}{\sin^3(\pi(\th-t))}|
\end{align*}
whence
\begin{align*}\D\tau(\th)&\le \|\D\mathfrak{b}\|_\infty+\pi
\frac{\sum_{t\in\mathfrak{s}_\phi}\s(\{t\})|\tfrac{\cos(\pi(\th-t))|}{\sin^3(\pi(\th-t))}|}{s'(\th)^2}\le
 \|\D\mathfrak{b}\|_\infty+\pi\sum_{t\in\mathfrak{s}_\phi}\frac{\tfrac{\s(\{t\})|\cos(\pi(\th-t))|}{|\sin^3(\pi(\th-t))|}}{
 (\tfrac{\s(\{t\})}{\sin^2(\pi(\th-t))})^2}\\ &=
 \|\D\mathfrak{b}\|_\infty+\tfrac{\pi}2\sum_{t\in\mathfrak{s}_\phi}\tfrac{|\sin(2\pi(\th-t))|}{\s(\{t\})}\le
 \|\D\mathfrak{b}\|_\infty+\tfrac{\pi}2\sum_{t\in\mathfrak{s}_\phi}\tfrac1{\s(\{t\})}\\ &=:M<\infty.\ \CheckedBox\ \text{\eqref{Adler}}
\end{align*}
\f{\bf Example}
\ \

We exhibit an inner function $\phi:\mathbb{D}\righttoleftarrow$ with $\phi(0)=0$,\ and singularity set $\mathfrak{s}_\phi=\{\chi(a),\chi(b)\}$ with $0\le a<b\le 1$; so that
$\tau=\tau(\phi):[a,b]\to [\tau(a+),\tau(b-)]\subset (0,1)$ is an homeomorphism. Thus the restriction is not an Adler arc map. Indeed, such an inner function must be a  Blaschke product since otherwise either  $\chi(a)$  is a singularity of the singular factor of $\phi$ and by \eqref{faMedkit} $\tau(a-)=0$, or $\chi(b)$ is and $\tau(b-)=1$.  This would contradict
$[\tau(a+),\tau(b-)]\subset (0,1)$.
\

To construct $\phi$ first define an inner function of $\Bbb R^{2+}$ (the upper half plane): $B:\Bbb R^{2+}\righttoleftarrow$ by
$$B(z):=\sum_{n=1}^\infty\tfrac1{2^{n+1}}(\tfrac{1+s_nz}{s_n-z}+\tfrac{1+t_nz}{t_n-z})$$ with $s_n:=-\tfrac1n$
and $t_n=1+\tfrac1n$; then $B$ is an inner function of $\Bbb R^{2+}$ and
$B(i)=i,\ B:(0,1)\to\Bbb R$ is continuous, increasing with $B(0+),\ B(1-)\in\Bbb R$. Let $\psi:\Bbb D\to\Bbb R^{2+},\ \psi(z)=i\tfrac{1+z}{1-z}$, then $\phi:=\psi^{-1}\circ B\circ\psi$ is as advertised with
$a=\tfrac12$ and $b=\tfrac34$.
\subsection{ Radial limit set}
\

The {\it radial limit set} of $\phi$ is
$$\Lambda_{\phi}:=\{\b\in\bdy\mathbb{D}:\ \phi(r\b)\xrightarrow[r\to 1-]{}\phi(\b)\in
\bdy\mathbb{D}\}.$$
For example,  if $\phi(z):=\exp[-\tfrac{1+z}{1-z}]$, then 
$$\phi(r\xi)\xrightarrow[r\to 1-]{}\begin{cases}& e^{-i\cot(\frac{\th}2)}\ \ \ \xi=\chi(\th)\ne 1;
\\ & 0\ \ \ \ \ \ \ \ \ \ \ \xi=1                                      
                                         \end{cases}$$

and $\Lambda_{\phi}=\bdy\mathbb{D}\setminus\{1\}$.

\subsubsection{ Angular derivatives}
\

The inner function $\phi:\mathbb{D}\righttoleftarrow$ has an
{\it angular derivative} $\b\in\Bbb C$ at $\xi\in\Lambda_{\phi}$ if
$$\exists\  \ \b\in\Bbb C\ \ \ \st\ \ \frac{\phi (z)-\phi(\xi)}{z-\xi}\xrightarrow[z\to\xi]{\angle}\b;$$ that is 
$$\frac{\phi (z)-\phi(\xi)}{z-\xi}\xrightarrow[z\to\xi,\ |\xi-z|\le K(1-|z|)]{}\b\ \ \forall\ K>0.$$
Denote the angular derivative at $\xi$ by $\b=:\phi_\angle'(\xi)$.

 \proclaim{Proposition 2.6\ \ \cite[\S3]{Saksman}}\ \ Let $\phi:\mathbb{D}\righttoleftarrow$ be inner
 with restriction  $\tau=\tau(\phi)$.\  The following are equivalent for $\xi=\chi(x)\in\bdy\mathbb{D}$
\sbul $\phi$ has an angular derivative at $\xi$;
\sbul $\phi'(z)\xrightarrow[z\to\xi]{\angle}\b\in\Bbb C\ \&$ in this case $\b=\phi_\angle'(\xi)$;
\sbul $\int_\Bbb T\tfrac{d\nu_w(t)}{|\xi-\chi(t)|^2}<\infty$ for some {\rm (hence all)}  $w\in\Bbb T,\ w\ne x$.\endproclaim
\proclaim{Proposition 2.7 \  \ {\rm\small\cite{Aleksandrov87}, also \cite[Theorem 9.6]{Saksman}.}}\label{AH} \

Let $\phi:\mathbb{D}\righttoleftarrow$ be inner and let $E\in\B(\Bbb T),\ m(E)>0$, then
$\tau=\tau(\phi)$ is locally invertible on $E$ iff $\phi$ has an angular derivative at $\chi(x)$
for $m$-a.e. $x\in E$.\endproclaim
\demo{Proof} \ \ $\Rightarrow$ Local invertibility on $E$ entails forward nonsingularity on $E$ whence, by \cite{Hei},  existence of angular derivatives a.s. on $E$.
 \
 
 To see  $\Leftarrow$, suppose that $\phi$ has an angular derivative at $\chi(x)$
for $m$-a.e. $x\in E$. By \cite[Lemma 1.5]{Craizer-E},  $\tau$ is {\it almost uniformly
differentiable} on $E$ in the sense that 
\sbul $\exists\ \ E_k\in\B(\Bbb T),\ E_k\uparrow\ E\ \mod m\ \st$:

$\forall\ k,\ \e>0\ \exists\ \d=\d(k,\e)>0$ so that
\begin{align*}\tag*{\dstechnical}\label{dstechnical}
 |\tau(x)-\tau(y)-(x-y)g(x)|\le \e|x-y|\ \forall\ x,y\in E_k,\ |x-y|<\d.
\end{align*}

where $g(x)=|\phi'(\chi(x))|$ with $\phi'$ the angular derivative of $\phi$ at $\chi(x)$.

\

To  see  that $\tau$ is forward nonsingular on $E$,  we note that by possibly shrinking the $E_k$ (as in \eqref{dstechnical}), we may assume  in addition that $\exists\ M_k>0$ ($k\ge 1$) so that $g\le M_k$ on $E_k$.
\

Suppose that $k\ge 1\ \&\ S\in\B(E_k),\ m(S)=0$ and fix $\e>0$.
There are intervals $\{I_n:\ n\ge 1\}$ so that
$$S\subset\bigcup_{n\ge 1}I_n,\ |I_n|\le\d(k,1)\ \& \ \sum_{n\ge 1}|I_n|<\tfrac{\e}{M_k+1}.$$
By \eqref{dstechnical},
$$|\tau(x)-\tau(y)|\le (g(x)+1)|x-y|\le (M_k+1)|x-y|$$
and $m(\tau(E_k\cap I_n))\le (M_k+1)m(I_n)$.
\

It follows that  $\tau(S)\subset\bigcup_{n\ge 1}\tau(E_k\cap I_n)$
whence 
$$m(\tau(S))\le \sum_{n\ge 1}m(\tau(E_k\cap I_n)\le (M_k+1)\sum_{n\ge 1}m(I_n)<\e$$
and $m(\tau(S))=0$.
\

For $S\in\B(E),\ m(S)=0$,
$$m(\tau(S))\xleftarrow[k\to\infty]{}m(\tau(E_k\cap S))=0$$  and
$\tau$ is forward nonsingular on $E$ whence locally invertible on $E$ by Theorem 2.2. \ \
\ \ \CheckedBox

An example of a probability  preserving restriction which is
"a.e. continuum to one" (a.e. Clark measure is nonatomic) was constructed in \cite{Donoghue} (see also \cite[Ex. 9.7]{Saksman}).

\

In particular, 
 for inner functions $\phi$, $\tau(\phi)$ is locally invertible on $\Bbb T$ iff $\phi$ has an angular derivative at a.e. point on $\bdy\mathbb{D}$.
 \
 
The next result  shows that this property is independent of the singularity set.

\proclaim{Proposition 2.8}\label{prop2.7}
\ 

Let $E\subseteq\Bbb T$ be a  closed set, then
\

$\exists$ an inner function $\phi:\mathbb{D}\righttoleftarrow$ with locally invertible restriction $\tau=\tau(\phi)$  so that $\phi(0)=0\ \&$
$\frak s_\phi=\chi(E).$\endproclaim
\demo{Proof}
\

Let $\G\subset\Bbb T$ be countable so that $\G'=E$.

\

We'll construct $\phi$ via a Clark measure.
\

First fix $\e:\G\to\Bbb R_+$ so that $\sum_{\g\in\G}\e(\g)<\infty$ and then fix $\pi\in\mathcal{P}(\G)$
so that
$$\sum_{\g\in\G}\tfrac{\pi_\g}{\e(\g)^2}<\infty.$$

 We claim that
\begin{align*}\tag*{\P1}\label{P1}\int_\Bbb T\tfrac{d\pi(t)}{|t-x|^2}<\infty\ \text{for $m$-a.e.}\ x\in\Bbb T
\end{align*}
where $\pi:=\sum_{\g\in\G}\pi_\g \d_\g\in\mathcal{P}(\mathbb{T})$.
\demo{Proof of \eqref{P1}}
\

Since 
$$\sum_{\g\in\G}m(B(\g,\e(\g)))=2\sum_{\g\in\G}\e(\g)<\infty$$
with $B(x,\e):=(x-\e,x+\e)$,
we have by the Borel-Cantelli lemma that

$\exists\ K\in\B(\Bbb T),\ K\cap\G=\emptyset,\ m(K)=1$ so that $\forall\ x\in K,\ 
\exists\ \G_0(x)\in\G$ finite, so that
\begin{align*}\tag*{\dsliterary}\label{dsliterary}
|x-\g|\ge\e(\g)\ \ \forall\ \g\notin \G_0(x).
\end{align*}
\

Let $x\in K$, then
\begin{align*}\int_\Bbb T\tfrac{d\pi(t)}{|t-x|^2}&=\sum_{\g\in\G}\tfrac{\pi_\g}{|x-\g|^2}
=(\sum_{\g\in\G_0(x)}+\sum_{\g\notin\G_0(x)})\tfrac{\pi_\g}{|x-\g|^2}. 
\end{align*}
Since $K\subset\Bbb T\setminus\G$, we have 
$$\sum_{\g\in\G_0(x)}\tfrac{\pi_\g}{|x-\g|^2}<\infty\ \forall\ x\in K$$ and
\begin{align*}\sum_{\g\notin\G_0(x)}\tfrac{\pi_\g}{|x-\g|^2}\le
\sum_{\g\notin\G_0(x)}\tfrac{\pi_\g}{\e(\g)^2}<\infty.\ \CheckedBox\ 
\P1 
\end{align*}
Next, we define $F:\mathbb{D}\to \Bbb C$ by
$$F(z):=\int_\Bbb T\tfrac{\chi(t)+z}{\chi(t)-z}d\pi(t).$$
Note that 
$$\text{\tt Re}(\tfrac{\chi(t)+z}{\chi(t)-z})=\tfrac{1-|z|^2}{|\chi(t)-z|^2}>0\ \forall\ z\in\mathbb{D}$$
so $F:\mathbb{D}\to\Bbb R_+\x\Bbb R$. Moreover, since $\pi(\Bbb T)=1$,
$F(0)=1$.

\

\Par2 \ To construct the inner function, define $\phi:=\tfrac{F-1}{F+1}:\mathbb{D}\to\Bbb C$.
\

Since, $F(0)=1,\ \phi(0)=0$ and since $\text{\tt Re}F>0$ on $\mathbb{D},\ \phi:\mathbb{D}\righttoleftarrow$.
\

Since $\pi\perp m$, we have that for $m$-a.e. $x\in\Bbb T$,
$$\exists\ \lim_{r\to 1-}F(r\chi(x))=:F(\chi(t))\in i\Bbb R$$
                                                            
whence
for such $x\in\Bbb T$

\

$$\phi(r\chi(x))\xrightarrow[r\to 1-]{}\tfrac{F(\chi(x))-1}{F(\chi(x))+1}\in\bdy\mathbb{D}$$
and $\phi:\mathbb{D}\righttoleftarrow$ is inner.

\

Moreover, $\nu_0=\pi$ whence by\ \ \eqref{P1}\ \  and Proposition 2.6, $\phi$ has an angular derivative at $\chi(t)$ for a.e. $t\in\Bbb T$. By Proposition 2.7, $\tau(\phi)$ is locally invertible.
 
\

To finish, we note that $\frak s_\phi=\text{\tt spt}\,\pi\ '=\G'=\chi(E)$.\ \CheckedBox

The Baire category situation is different:
\proclaim{ Proposition 2.9}
\

\ Suppose that the inner function $\phi:\mathbb{D}\righttoleftarrow$ has an angular derivative at a residual set of points,
then
$\frak s_\phi$ is nowhere dense.\endproclaim\demo{Proof}\ \ 
\

Fix  $w\in\Bbb T$  and define $F:\Bbb T\to (0,\infty]$ by
$$F(x):=\int_\Bbb T\tfrac{d\nu_w(t)}{|\chi(x)-\chi(t)|^2}$$
and, for  $r\in (0,1)$,\ define $F_r:\Bbb T\to\Bbb R_+$ by
$$F_r(x):=\int_\Bbb T\tfrac{d\nu_w(t)}{|r\chi(x)-\chi(t)|^2},$$
then for each  $0<r<1$,\  $F_r:\Bbb T\to\Bbb R_+$ is continuous  and \ \ 
\f$F_r(x)\xrightarrow[r\to1-]{}F(x)\in (0,\infty]\ \forall\ x\in\Bbb T$ by dominated convergence in case $F(x)<\infty$ and by Fatou's lemma in case $F(x)=\infty$

\

By Baire's simple limit theorem,  
$$C_\infty:=\{x\in\Bbb T:\ F:\Bbb T\to (0,\infty]\ \text{ is continuous at}\ x\}$$ is residual in $\Bbb T$.

By assumption, 
$\Lambda:=[F<\infty]$
is residual, whence so is $C:=C_\infty\cap\Lambda$.

Thus $\forall\ x\in C,\ \exists\ 0<a_x<b_x$ and an open interval $J_x\ni x$ so that 
$$a_x<F<b_x\ \text{on}\ J_x.$$

By  Seidel's theorem   \cite{Seidel}
(see also \cite[Thm. 7.48]{Zygmund}), 
\

$$\frak s_\phi\cap\chi(J_x)=\emptyset,$$ whence $\frak r_\phi$ is open and dense and 
$\frak s_\phi$ is nowhere dense. 
\ \CheckedBox

\

\

\

\section{ Quasicompactness of transfer operators}\label{QcTO}
\

\

Let $\mathcal{L}$ be a Banach space. An operator  $P\in\hom(\mathcal L,\mathcal L)$  is called {\it quasicompact} if
$\exists$ $A=A(P)\in\hom(\mathcal L,\mathcal L)$ of form
$$A=\sum_{k=1}^N\l_kE_k$$ with $N\ge 1,\ E_1,\dots,E_N\in\hom(\mathcal L,\mathcal L)$ finite dimensional  projections,
$\l_1,\dots,\l_N\in \mathbb{S}^1:=\{z\in\Bbb C:\ |z|=1\}$
so that  the spectral radius
$$\rho(P-A):=\lim_{n\to\infty}\|(P-A)^n\|^{\frac1n}_{\hom(\mathcal L,\mathcal L)}<1.$$
\

Let $(X,m,T)$ be a nonsingular transformation with transfer operator $\T:L^1(m)\righttoleftarrow$.
We look for Banach spaces $\mathcal{L}\subset L^1(m)$ on which $\T:\mathcal{L}\righttoleftarrow$ acts  quasicompactly.

In case $(X,m,T)$ is a weakly mixing, probability preserving transformation
with  transfer operator $\T$ acting  quasicompactly. on $\mathcal{L}\subset L^1(m)$, then
$A(\T)f=\Bbb E(f)$ and  $\forall\ \th\in (\rho(\T-\Bbb E),1),\ \exists\ M>0$ so that
\begin{align*}\tag*{\dsmathematical}\label{dsmathematical}
 \|\T^n f-\Bbb E(f)\|_{\mathcal{L}}\le M\th^n\|f\|_{\mathcal{L}}\ \forall\ f\in\mathcal{L}.
\end{align*}
The property \eqref{dsmathematical} is aka {\it exponential decay of correlations} as it entails
$$|\int_X u\cdot v\circ T^ndm-\Bbb E(u)\Bbb E(V)|\le M\th^n\|u\|_{\mathcal{L}}\|v\|_{L^1(m)}.$$
\

\subsection{Doeblin-Fortet operators on an adapted pair}\label{quasicompact}
\

Let   $\mathcal L\subset \mathcal{C}\subset L^1(m)$ be  Banach spaces  so that
                                                                               
$(\mathcal{C},\mathcal{L})$ is an 
{\it adapted pair} in the sense that
\sbul 
$\|\cdot\|_{{L^1(m)}}\le\|\cdot\|_{\mathcal C}\le\|\cdot\|_{\mathcal L},\ (\o {\mathcal L})_{{L^1(m)}}={L^1(m)}$,
and $\mathcal L$-closed, bounded sets are $\mathcal{C}$-compact.
\

For example both  $(L^1(m),\text{\tt Lip}(\Bbb T))$  and  $(L^1(m),\text{\tt BV}(\Bbb T))$  are adapted pairs where
\sbul $\text{\tt Lip}(\Bbb T)$ denotes the Lipschitz functions on $\Bbb T$
(equivalently the absolutely continuous functions with essentially bounded derivative); with norm $\|f\|_{\text{\tt \tiny Lip}}:=\|f\|_1+\|f'\|_\infty$; {\tt Lip}-closed, bounded sets being $L^1$-compact by the Arzela-Ascoli theorem; and
\sbul $\text{\tt BV}(\Bbb T)$ denotes the  functions of bounded variation on $\Bbb T$

with norm $\|f\|_{\text{\tt \tiny BV}}:=\|f\|_1+\bigvee f$ where
$$\bigvee f:-\sup\,\{\sum_{k=0}^{n-1}|f(t_{k+1})-f(t_k)|:\ 0<t_1<t_2<\dots<t_n=1\};$$ {\tt BV}-closed, boiunded sets being
$L^1$-precompact by Helly's theorem.
\

 As in \cite[Chapter 3]{Norman}, we say that an operator $P\in\hom(\mathcal L,\mathcal L)\cap\hom(\mathcal{C},\mathcal{C}) $ is {\it Doeblin-Fortet} ({\tt D-F}) on $(\mathcal{C},\mathcal{L})$  if
 \begin{align*}&\tag*{{\tt DF}(i)}\label{DF(i)}
 \|P^nf\|_\mathcal{C}\le H\|f\|_\mathcal{C}\ \forall n\in\Bbb N,\ f\in L^1(m)\ \&\\ &\tag*{{\tt DF}(ii)}\label{DF(ii)}\ \exists\ \kappa\ge 1\ \st\
\|P^\kappa f\|_{\mathcal L}\le \th \|f\|_{\mathcal L}+R\|f\|_\mathcal{C}\ \forall\ f\in \mathcal L.  
 \end{align*}

where $R,H\in\Bbb R_+$ and $\th\in (0,1)$.

\

\f{\bf Example 3.2: Adler arc maps}\label{prop3.2}
\

Let $(\Bbb T,m=\text{\tt\small Leb},T)$ be an Adler map (as on p. \pageref{Exp}). It is a well known follore result that
  the transfer operator $\T$  is {\tt D-F} on the adapted pair   $(L^1(m),\text{\tt Lip}(\Bbb T))$.

It is also {\tt D-F} on   $(L^1(m),\text{\tt BV}(\Bbb T))$
because, by \cite[Corollary 1]{Roland} an Adler map satisfies the assumptions of \cite[Proposition 1]{Rychlik}
which proves the {\tt D-F} inequality on  $(L^1(m),\text{\tt BV})$.

\

The following lemma is a well-known consequence of the Yosida-Kakutani mean ergodic theorem
(\cite[Theorem1]{Yosida-Kakutani}). See also
\cite{IT-M}, \cite[Chapter 3]{Norman}, \cite{HH, Par-Pol}, \cite[Thm 1]{L-Y}.
\proclaim{Lemma 3.1}
\

Suppose that $P$ is a Doeblin-Fortet operator on the adapted pair $(\mathcal{C},\mathcal{L})$.
\

If $f\in\mathcal{C},\ Pf=f$, then $f\in\mathcal{L}$.\endproclaim


It is shown in \cite{IT-M} (see also
\cite[Chapter 3]{Norman}, \cite{HH, Par-Pol}) that a D-F operator $P\in\hom(\mathcal L,\mathcal L)$  has spectral radius $\rho(P)\le 1$ and that,
if $\rho(P)=1$, then $P$ is  quasicompact.
\

\subsubsection{Quasicompactness 
and the Central Limit Theorem}
\

If $\psi\in\mathcal{L},\ \Bbb E(\psi)=0$, then by Leonov's theorem (\cite{Leo})
\begin{align*}\tag{{\tt Leonov}}\label{Leonov}
 \begin{split}
  &\exists\ \lim_{n\to\infty}\tfrac1n\Bbb E(\psi_n^2)=:\s^2_{\psi}\ge 0\ \text{where}\ \psi_n:=\sum_{k=0}^{n-1}\psi\circ T^k;\\ &
  \text{with equality iff}\ \psi=g-g\circ T\ \text{for some}\ g\in\psi.
 \end{split}
\end{align*}
If, in addition, $\s_\psi>0$, then (\cite{Gordin-Martingale}) the {\it stationary process} $(X,m,T,\psi)$
satisfies the central limit theorem:
\begin{align*}\tag*{{\tt CLT}}\label{CLT}
m([\frac{\psi_n}{\s_\psi\sqrt n}\le t])\xrightarrow[n\to\infty]{}\tfrac1{\sqrt{2\pi}}\int_{-\infty}^te^{-\frac{s^2}2}ds.
\end{align*}

\subsection{Hardy spaces}\label{hardy}
\

The {\it harmonization} of  $f\in L^p(m)$ ($1\le p\le\infty$) is $\widetilde{f}:\mathbb{D}\to\Bbb C$ is defined by
$$\widetilde{f}(z):=\int_\Bbb T p_zfdm$$ where $p_z$ is as in \eqref{PHtunny} on p. \pageref{PHtunny}.
\

It is harmonic in $\mathbb{D}$ and satisfies
$$\sup_{r\in (0,1)}\|\widetilde{f}(r\chi)\|_p=\|f\|_p.$$
It is classical that the Hardy spaces consist of harmonizations:
\begin{align*}&h^p(\mathbb{D}):=\{F:\mathbb{D}\to\Bbb C\ \text{\small harmonic},\ \sup_{r\in (0,1)}\|F(r\chi)\|_p<\infty\}=\{\widetilde{f}:f\in L^p(m)\};
\\ &H^p(\mathbb{D}):=
\{f\in h^p:\ \widetilde{f}\ \text{\small  analytic on}\ \mathbb{D}\}\cong\{f\in L^p(\Bbb T):\ \widehat{f}(n)=0\ \forall\ n<0\}.\end{align*}

   Let 
   $$\Lambda_{\widetilde{f}}:=[\exists\ \lim_{r\to 1-}\widetilde{f}(r\chi)=:\widetilde{f}(\chi)\in\Bbb C],$$ then by Fatou's theorem, $m(\Lambda_{\widetilde{f}})=1\ 
   \&\ \widetilde{f}(\chi)=f$ a.e..

\subsubsection{Action of the transfer operator}
\

Let $\phi:\mathbb{D}\righttoleftarrow$ be inner  with restriction $\tau=\tau(\phi )$, then
(\cite{Aleksandrov87}, see also \cite[Theorem 3.1]{Saksman}):
$\ttau H^p(\mathbb{D})\subset H^p(\mathbb{D})$ and  if $\phi (0)=0$, then for $d\ge 1,\ \ttau(\chi^d)$ is a polynomial in $\chi$ of degree at most $d$. Moreover
\proclaim{Lemma  3.3}
\

Let $\phi:\mathbb{D}\righttoleftarrow$ be inner with  $\phi^{(k)}(0)=0\ \forall\ 0\le k<\kappa$,  then for $N\ge 1$:
$\ttau^N(\chi^d)=0\ \forall\ 1\le d<\kappa^N$ and for
$d\ge \kappa^N,\ \ttau^N(\chi^d)=\sum_{\ell=1}^{\lfl\frac{d}{\kappa^N}\rfl}a^{(N}(d,\ell)\chi^\ell$ where
$a^{(N)}(k,\ell)=\overline{\widehat{(\phi^{[N]})^\ell}}(k)$. \endproclaim

\demo{Proof}\ \
\

Let $d\ge 1$ and write $\ttau^N(\chi^d)=\sum_{\ell\in\Bbb Z}a^{(N)}(d,\ell)\chi^\ell$, then
\begin{align*}a^{(N)}(d,\ell)&=\< \ttau^N(\chi^d),\chi^\ell\>=
\<\chi^d,\chi^\ell\circ T^N\>=
\<\chi^d,\phi^{N}(\chi)^\ell\>=\overline{\widehat{(\phi^N)^\ell}(d)}.
\end{align*}
\footnote{\tiny Note that here $\phi^{[N]}:=\underbrace{\phi\circ\dots\circ \phi}_{\text{\tiny $N$ times}}$ whereas
$(\phi^{[N]})^\ell:=\underbrace{\phi^{[N]}\cdot\dots\cdot \phi^{[N]}}_{\text{\tiny $\ell$ times}}$.}
\

Thus
\begin{equation*}\tag*{\dsaeronautical}\label{dsaeronautical}  a^{(N)}(d,\ell)=0\ \text{unless}\ \ell\ge 1\ \&\ d\ge \kappa^N\ell.\ \ \CheckedBox
\end{equation*}

\subsection{Weighted  Hilbert spaces}\label{WHS}
\

A (Hilbert space) {\it weight} is a sequence   $\u w\in\Bbb R_+^{\Bbb N_0}$, define satisfying $1=w(0)<w(1)<\dots <w(n)\uparrow\infty$. The associated
{\it weighted Hilbert space} is
$$\mathcal{h}_{\u w}:=\{f\in L^2(m):\ \|f\|_{\u w}^2:=\sum_{n\in\Bbb Z}w(|n|)|\widehat{f}(n)|^2<\infty\},$$
 equipped with the inner product
$$\<u,v\>_{\u w}:=\sum_{n\in\Bbb Z}w(|n|)\widehat{u}(n)\overline{\widehat{v}(n)}.$$

\f{\bf Classical Examples}
\

\sms (i) For $w(n)=n^2,\  \mathcal{h}_{\u w}$ is isomorphic to the Sobolev space:
$$W^{1,2}(\Bbb T):=\{f\in C(\Bbb T):\ f\ \text{\small a.c.}\ \&\ f'\in L^2(m)\}.$$
\sms (ii)\ For $b>1$, let $w_{1,b}(n):=b^n$, then
$$\mathcal{k}_b=\mathcal{h}_{\u w_{1,b}}\cong h^2(B_\Bbb C(0,b)^o),$$
\

We'll call a weight $\u w$ {\it summable} if $\sum_{n\ge 1}\tfrac1{w(n)}<\infty$. Both examples above are summable.
\

\proclaim{Proposition 3.4}
\

\ \ If
$\u w$ is a summable weight, then
 $(L^2(\Bbb T),\mathcal{h}_{\u w})$ is an adapted pair.\endproclaim
\demo{Proof} \ \ We show that $B(R):=\{f\in \mathcal{h}_{\u w}:\ \|f\|_{\u w}\le R\}$ is strongly compact in $L^2(m)$.
To see this let $f_j\in B(R)$ ($j\ge 1$), then for $n\in\Bbb Z,\ j\ge 1$,
$$|\widehat{f}_j(n)|\le\tfrac{\|f_j\|_{\u w}}{\sqrt{w(|n|)}}\le \tfrac{R}{\sqrt{w(|n|)}}$$ and $\exists\ j_\ell\to\infty,\ a\in\ell^2(\Bbb Z)$ so that
$$\widehat{f}_{j_\ell}(n)\xrightarrow[\ell\to\infty]{}a(n).$$
We claim that $\sum_{n\in\Bbb Z}w(|n|)|a(n)|^2\le R^2$.
To see this,
\begin{align*}R^2&\ \ge \sum_{n\in\Bbb Z}w(|n|)|\widehat{f}_{j_\ell}(n)|^2\ \ge\  \sum_{|n|\le N}w(|n|)|\widehat{f}_{j_\ell}(n)|^2\ \ \forall\ N\ge 1,\\ &
\xrightarrow[\ell\to\infty]{}\sum_{|n|\le N}w(|n|)||a(n)|^2\ \xrightarrow[N\to\infty]{}\sum_{n\in\Bbb Z}w(|n|)||a(n)|^2.
\end{align*}
Let $A:=\sum_{n\in\Bbb Z}a(n)\chi^n\in \mathcal{h}_{\u w}$ with $\|A\|_{\u w}\le R$.
\

To see that $f_{\ell}\xrightarrow[\ell\to\infty]{L^2(m)}\ A$, by  the Riesz-Fischer theorem,
\begin{align*}\|f_{j_\ell}-A&\|_{L^2(m)}^2=
\sum_{n\in\Bbb Z}|\widehat{f}_{j_\ell}(n)-a(n)|^2\  \xrightarrow[\ell\to\infty]{}0\\ &\text{\Large$\mathbf{\because}$}\ \ 0\xleftarrow[\ell\to\infty]{}|\widehat{f}_{j_\ell}(n)-a(n)|^2\ll \tfrac1{w(|n|)}.\ \ \CheckedBox
\end{align*}
\

The rest of this section is devoted to showing that the transfer operators of probability preserving,
non-M\"obius, inner functions act quasicompactly on certain weighted Hilbert spaces.

 Ivrii and Urbanski (\cite{OlegMar}) obtained {\it i.a.} spectral gaps for  the action of $\ttau$ on $W^{1,2}(\Bbb T)$ ($\tau=\tau(\phi)$ with $\phi$ inner, $\phi(0)=0$)  and we obtain {\it i.a.} them on $\mathcal{k}_b$ (Proposition 3.5 below). In both cases,
 the minimal essential radius (as in \S\ref{ess-specrad}) is the ''Koenigs eigenvalue`` $|\phi'(0)|$
 (Proposition 3.7 below). However we obtain superexponential decay of correlations on (e.g. $\mathcal{k}_b$) when $\phi'(0)=0$ (Proposition 3.6 below).

\subsubsection{Admissible weighted Hilbert spaces}\
\

Call  a summable weight $\u w\in\Bbb R_+^\Bbb N$ and its associated
weighted Hilbert space $\mathcal{h}_{\u w}$
  {\it admissible} if

\begin{align*}\tag*{\faLinux}\label{faLinux}
\exists\ C=C_{\u w}>0\ \text{s.t.}\  W(\lfl\tfrac{n}K\rfl)\le C\, \tfrac{w(n)}{w(K)}\ \forall\ n\ge K\ge 1\end{align*}
where $W(n):=\sum_{k=1}^nw(k)$.

We'll call any $C_{\u w}$ satsfying \eqref{faLinux} an {\it admissibility constant } for $\u w$.
\

For example, for $b>1$, $w_{1,b}(n)=b^n$ defines an admissible weight with e.g. $C_{\u w_{1,b}}=\tfrac{b}{b-1}$.
\
Also, for $b>1,\ s> 1$,   $\u w_{s,b}$ defined by $w_{s,b}(n):=b^{n^s}$ also defines an admissible weight.

On the other hand,  for $t>0,\ \u v_t$ defined by $v_t(n):=n^t$ is not admissible
(although summable for $t>1$).
\

Recall from \cite[Definition I.2.10]{Izzy} that a Banach space
$B\subset L^1(\Bbb T,m)$ is {\it homogeneous} if
\begin{align*}\tag*{\dschemical}\label{dschemical}
\begin{split}
&f\in B,\ s\in\Bbb T\ \implies\ f_s\in B,\
 \ \|f_s\|_B=\|f\|\\ & \ \&\ \|f-f_s\|_B\xrightarrow[s\to 0]{}0\  \text{with}\ f_s(x):=f(x-s);.\end{split}
\end{align*}
Consequently (\cite[Theorem I.2.11]{Izzy}), if $B$ is homogeneous, then for $f\in B$,
\begin{align*}\tag*{\dsagricultural}\label{dsagricultural}
f*p_r\in B,\  \|f*p_r\|_B=\|f\|_B\ \forall\ 0<r<1\ \&\ f*p_r\xrightarrow[r\to 1]{B}f.
\end{align*}
Any summably weighted Hilbert space $\mathcal{h}_{\u w}$ is  homogeneous.

\


\proclaim{ Proposition 3.5  (exponential decay of correlations)}
\

Let $\phi:\mathbb{D}\righttoleftarrow$ be inner, non-M\"obius with $\phi(0)=0\ \&\ \tau=\tau(\phi)$ and let
$\u w\in\Bbb R_+^\Bbb N$ be admissible.
\

If $\exists\ 0<R=R_{\u w}<1$ so that $\sum_{n\ge 1}\tfrac1{R^{2n}w(n)}<\infty$, then for $\forall\ |\phi'(0)|<\rho<1,\ \exists\ M>0$ so that
\begin{align*}\tag*{\Wheelchair}\label{Wheelchair}
\|\ttau^Nu-\Bbb E(u)\|_{\u w}\le M\rho^N\|u\|_{\u w}\ \forall\ u\in \mathcal{h}_{\u w},\ N\ge 0
\end{align*} where $\Bbb E(u):=\int_\Bbb Tudm=\widehat{u}(0)$.

\

\endproclaim\demo{Proof}
\

Fix $\rho\in (|\phi'(0)|,1)$.  We first show that
\begin{align*}\tag*{\phone}\label{phone}\forall\ 0<r<1,\ \exists\ M>0\ \text{s.t.}\
|\phi^{[N]}(z)|\le M\rho^N\ \forall\ z\in\mathbb{D},\ |z|\le r
\end{align*}
where $\phi^{[N]}=\underbrace{\phi\circ\phi\circ\dots\circ \phi}_{\text{\tiny $N$ times}}$.
\demo{Proof of (\eqref{phone})}\ \  By assumption $\phi(z)=zg(z)$ where $g:\mathbb{D}\righttoleftarrow$ is inner. Thus $\phi'(0)=g(0)$ and $\exists\ \frak r=\frak r_\rho$ so that
$$|g(z)|\le\rho\ \&\ |\phi(z)|\le\rho |z| \forall\ |z|\le\frak r.$$
\

Next, $\exists\ N_\rho$ so that $|\phi^{[N_\rho]}(z)|\le\frak r\ \forall\ |z|\le r$ whence
for $|z|\le r$
$$|\phi^{[n]}(z)|\le \rho^{n-N_\rho}\frak r\ \forall\ n>N_\rho$$ and \eqref{phone} follows.\ \CheckedBox

\

\demo{Proof of \eqref{Wheelchair} for $u\in\mathcal{h}_{\u w}\cap H_0^2$}
\

For $u\in\mathcal{h}_{\u w}\cap H^2_0$ set $v_n:=\tfrac{\widehat{u}(n)}{R^n}$ where $R=R_{\u w}$, then
$$\sum_{n\ge 1}|v_n|=\sum_{n\ge 1}\sqrt{w(n)}
|\widehat{u}(n)|\cdot\tfrac1{R^{n}\sqrt{w(n)}}\le \|u\|_{\u w}^2\sum_{n\ge 1}\tfrac1{R^{2n}w(n)}<\infty,$$
whence $v:=\sum_{n\ge 1}v_n\chi^n\in C(\Bbb T)$ and $v*p_R=u$.

Now $$u=v*p_R=\int_\Bbb Tv(t)p_{R\chi(t)}dm$$ whence
\begin{align*}u_N:=\ttau^N u&=\int_\Bbb Tv(t)\ttau^N(p_{R\chi(t)})dm\\ &=
 \int_\Bbb Tv(t)p_{\phi^{[N]}(R\chi(t))}dm\ \text{by \eqref{dsheraldical}}
\end{align*}
and 
\begin{align*}\widehat{u}_N(\ell)&= \int_\Bbb T\int_\Bbb Tv(t)p_{\phi^{[N]}(R\chi(t))}\chi^{-\ell}dmdm\\ &=
1_\Bbb N(\ell)\int_\Bbb Tv(t)\widehat{p}_{\phi^{[N]}(R\chi(t))}(\ell)dm\\ &=
1_\Bbb N(\ell)\int_\Bbb Tv(t)(\overline{\phi^{[N]}(R\chi(t))})^\ell dm.
\end{align*}
      
      Thus
 \begin{align*}|\widehat{u}_N(\ell)|^2&\le 1_{\Bbb N}(\ell)\int_\Bbb T|v|^2dm
 \int_\Bbb T|\phi^{[N]}(R\chi(t))|^{2\ell}dt  \\ &\le
 \|u\|_{\u w}^2(M\rho^N)^{2\ell}.
 \end{align*}
 where $M>0$ is as in \eqref{phone}.
 \
 
   Let $N_0$ be so that $\frak b_n:=M^2\rho^{2n}<\tfrac12$ for $n\ge N_0$.
   For $N\ge N_0$
   \begin{align*}\|\ttau^N u\|^2_b &=\sum_{\ell\ge 1}b^\ell|\widehat{u_N}(\ell)|^2\le
   \|u\|^2_b\sum_{\ell\ge 1}\frak b_N^\ell\\ &=
   \tfrac{\frak b_N}{1-\frak b_N}\|u\|^2_b\le
   2M^2b\rho^{2N}\|u\|^2_b. \ \ \CheckedBox \  \text{\eqref{Wheelchair} }
   \end{align*}

  To continue let $w\in \mathcal{h}_{\u w}$, then
  $$w=F+\overline{G}+\Bbb E(w)\ \text{for some}\ F,\ G\in\mathcal{h}_{\u w}\cap H_0^2$$
  and
  $$\ttau^N w=\ttau^N F+\overline{\ttau^N G}+\Bbb E(w).$$
By \eqref{Wheelchair} for $F\ \&\ G$,
  \begin{align*}\|\ttau^Nw-\Bbb E(w)\|^2_b&=\|\ttau^N F\|^2_b+\|\ttau^N G\|^2_b\\ &\le
   M\rho^{2N}(\|F\|^2_b+\|G\|^2_b)\\ &\le  M\rho^{2N}\|w\|_{\u w}^2.\ \ \text{\CheckedBox\ \eqref{Wheelchair}}
  \end{align*}

 \
 \proclaim{Proposition 3.6 (superexponential  decay of correlations)}

 Suppose that $\phi(z)=z^\kappa\Phi(z)$ with $\kappa > 1$ and $\Phi:\mathbb{D}\righttoleftarrow$ inner and let
$\u w\in\Bbb R_+^\Bbb N$ be admissible, then
\begin{align*}\tag*{\Bicycle}\label{Bicycle}
\|\ttau^Nu-\Bbb E(u)\|_{\u w}\le
\tfrac{\sqrt{C_{\u w}}}{\sqrt{w(\kappa^N)}}\|u\|_{\u w}\ \forall\ u\in \mathcal{h}_{\u w},\ N\ge 0.
\end{align*}\endproclaim

\demo{Proof}
\

By the Lemma 3.3,  for $N,\ k\ge 1$
$$\ttau^N(\chi^k)=\sum_{\ell=1}^{\lfl\frac{d}{\kappa^N}\rfl}a^{(N)}(k,\ell)\chi^\ell$$ where
$a^{(N)}(k,\ell)=\overline{\widehat{(\phi^N)^\ell}(k)}$ and $\sum_{\ell=1}^0=0$.

Thus also $\sum_{k\ge\kappa^N\ell}|a^{(N)}(k,\ell)|^2= 1.$

\demo{\P1: Proof of \eqref{Bicycle} for $u\in \mathcal{h}_{\u w}\cap H_0^2$}
\begin{align*}\ttau^N u&=\sum_{k\ge 1}u_k\ttau^N(\chi^k) \ \text{where}\ u_k=\widehat{u}(k),
\\ &=\sum_{\ell\ge 1,\ k\ge\kappa^N\ell}u_ka^{(N)}(k,\ell)\chi^\ell
\\ &=\sum_{\ell\ge 1}(\sum_{k\ge\kappa^N\ell}u_ka^{(N)}(k,\ell))\chi^\ell
 \end{align*}
 Thus,  using Cauchy-Schwartz and (\eqref{dsaeronautical}),
 \begin{align*}
   |(\ttau^Nu)_\ell|^2&=|\sum_{k\ge\kappa^N\ell}u_ka^{(N)}(k,\ell)|^2\\ &\le
  \sum_{k\ge \kappa^N\ell} |u_k|^2
 \end{align*}
and
  \begin{align*}\|\ttau^N u\|_{\u w}^2 &=\sum_{\ell\ge 1}w(\ell) |(\ttau^Nu)_\ell|^2\ \le
  \sum_{\ell\ge 1}w(\ell)   \sum_{k\ge \kappa^N\ell} |u_k|^2
  \\ &\le
  \sum_{k\ge \kappa^N}|u_k|^2\sum_{1\le \ell\le \tfrac{k}{\kappa^N}}w(\ell)\ =\sum_{k\ge\kappa^N}|u_k|^2W(\lfl\tfrac{k}{\kappa^N}\rfl)\\ & =
 \sum_{k\ge\kappa^N}|u_k|^2w(k)\cdot \tfrac{W(\lfl\tfrac{k}{\kappa^N}\rfl)}{w(k)}\\ &\le \tfrac{C_{\u w}}{w(\kappa^N)}\|u\|_{\u w}^2\ \text{by \eqref{faLinux}}.\ \ \CheckedBox\ \text{\P1}
  \end{align*}
  To continue let $u\in \mathcal{h}_{\u w}$, then
  $$u=F+\overline{G}+\Bbb E(u)\ \text{for some}\ F,\ G\in \mathcal{h}_{\u w}\cap H_0^2$$
  and
  $$\ttau^N u=\ttau^N F+\overline{\ttau^N G}+\Bbb E(u).$$
By \P1,
  \begin{align*}\|\ttau^Nu-\Bbb E(u)\|^2_{\u w}&=\|T^N F\|^2_{\u w}+\|\ttau^N G\|^2_{\u w}\\ &\le
   \tfrac{C_{\u w}}{w(\kappa^N)}(\|F\|^2_{\u w}+\|G\|^2_{\u w})\\ &=  \tfrac{C_{\u w}}{w(\kappa^N)}\|u\|_{\u w}^2.\ \ \text{\CheckedBox\ \eqref{Bicycle}}
  \end{align*}

\subsubsection{Essential spectral radius}\label{ess-specrad}
\

For $(\Bbb T,m,T)$ a weakly mixing, probability preserving transformation and suppose that $\T$ is a Doeblin-Fortet operator on the adapted pair $(L^1(\Bbb T,m),\mathcal{L})$.

 The {\it essential spectral radius} of $\T:\mathcal{L}\righttoleftarrow$ is
$$\rho_{\tt\tiny ess.}(\T,\mathcal{L}):=\rho(T,\mathcal{L}_0)=\lim_{n\to\infty}\|\T^n\|_{\hom(\mathcal{L}_0,\mathcal{L}_0)}^{\frac1n}$$
  where $\mathcal{L}_0:=\{f\in \mathcal{L}:\ \Bbb E(f)=0\}$.
  Equivalently, $\rho(\T,\mathcal{L})$ is the greatest lower bound of the collection of $\th\in (0,1)$ satisfying \eqref{dsmathematical}.

\proclaim{Proposition 3.7  \ (minimal essential spectral radius)}\ \footnote{c.f. \cite{BCJ}}
\

Let $\phi:\mathbb{D}\righttoleftarrow$ non-M\"obius, inner with $\phi(0)=0$, and suppose that $\T$ is a Doeblin-Fortet operator on the adapted pair
$(L^1(m),\mathcal{L})$ where $\exists\ R>0$ so that
$p_z\in\mathcal{L}\ \forall\ z\in\Bbb D,\ |z|<R$, then $\rho_{\tt\tiny ess.}(\T,\mathcal{L})\ge |\phi'(0)|$.
\endproclaim
\demo{Proof}
\

By
Koenigs' theorem, (\cite{Koenigs}, see also \cite[\S6.1]{Shapiro}),\ $\exists\ z\in\mathbb{D},\ |z|<R$   s.t.
$|\phi^n(z)|^{\frac1n}\xrightarrow[n\to\infty]{}|\phi'(0)|\ \&$ by \eqref{dsheraldical} (p.\pageref{dsheraldical}), if $\ttau$ satisfies
\eqref{dsmathematical} (p.\pageref{dsmathematical}) with constant $\th$ on $\mathcal{L}\in\mathfrak{h}$, then
$$M\th^n\ge \|\ttau^n(p_z)-1\|_\mathcal{L}\ge\|\ttau^n(p_z)-1\|_1\ge\e |\phi^n(z)|=|\phi'(0)|^{n+o(n)}.\ \CheckedBox$$
\section{ Perturbations and central limits}\label{pert}
\

Let $(X,m,T)$ be a weakly mixing, probability preserving transformation and let $\psi:X\to\Bbb C$ be measurable.
\

For $z\in\Bbb C$ so that $e^{z\psi}\in L^\infty(m)$, define the {\it perturbed {\rm (aka twisted)} transfer operator}
$\Pi_{z,\psi}:L^1(m)\righttoleftarrow$ by
$$\Pi_{z,\psi}f:=\T(e^{z\psi}f).$$
\

The operators $P_t=P_{t,\psi}:=\T(e^{it\psi}f)$ aka
{\it characteristic function operators} as $\Bbb E(P_t\mathbb{1})=\Bbb E(e^{it\psi})$.

Let $\mathcal{L}$ be a Banach space of functions on $\Bbb T$. A function $f:\Bbb T\to\Bbb R$ is an $\mathcal{L}$-{\it multiplier}
if $f\cdot u\in\mathcal{L}\ \forall\ u\in\mathcal{L}$. Let $M(\mathcal{L}):=\{\mathcal{L}-\text{\rm multipliers}\}$. By   the Resonance Theorem
$$\|f\|_{M(\mathcal{L})}:=\sup\,\{\|fu\|_{\mathcal{L}}:\ u\in\mathcal{L},\ \|u\|_{\mathcal{L}}=1\}<\infty\ \forall\ f\in M(\mathcal{L}).$$

Evidently, if  $\mathbbm{1}\in\mathcal{L}$, then $M(\mathcal{L})\subseteq\mathcal{L}$. Indeed, for $\mathcal{L}=${\tt Lip} or {\tt BV}, $M(\mathcal{L})=\mathcal{L}$.
\

In general, if $f\in M(\mathcal{L})$, then $f^N\in M(\mathcal{L})\ \forall\ N\ge 1$ with $\|f^N\|_{M(\mathcal{L})}
\le \|f\|_{M(\mathcal{L})}^N$ and $e^{zf}\in M(\mathcal{L})\ \forall\ z\in\Bbb C$ with $\|e^{zf}\|_{M(\mathcal{L})}\le
e^{|z|\|f\|_{M(\mathcal{L})}}$.
\

\proclaim{Theorem 4.1 (Nagaev's Theorem {\rm\small \cite{Nag,Rou})})}\label{Nagaev}\
\footnote{see also \cite{Par-Pol,HH}, \cite[Lemma 4.2]{A-De2}}
\

Suppose that $\T$ is a Doeblin-Fortet operator on the adapted pair
$(L^1(m),\mathcal{L})$ and that $\psi\in M(\mathcal{L})\cap\mathcal{L}$ satisfies $\Bbb E(\psi)=0$ and
\begin{align*}\tag*{$\clubsuit$}\label{clubsuit}
\begin{split}
&\exists\ \e>0\ \text{so that}\ P_{t,\psi}\in\hom(\mathcal{L},\mathcal{L}) \forall\ |t|<\e;\ \&\\ & t\mapsto P_{t,\psi}\ \text{is}\ C^2\ \ (-\e,\e)\to\hom(\mathcal{L},\mathcal{L}),
\end{split}
\end{align*}

\

then
\f {\rm (i)}\ \  $\exists\ 0<\mathcal{E}<\e$ so that $P_t$ is a Doeblin-Fortet operator on

$(L^2(m),\mathcal{L})\ \forall\ |t|<\mathcal{E}$.
\f{\rm (ii)}\ \ There are constants $K>0$ and $\th\in (0,1)$; and  $C^2$ functions
$\l:B(0,\mathcal{E})\to B_{\Bbb C}(0,1),\  N:B(0,\mathcal{E})\to\hom(\mathcal{L},\mathcal{L})$ so that
\begin{align*}\tag*{\faHome}\label{faHome}
 \|P_t^nh-\l(t)^n N(t)h\|_{\mathcal{L}}\le K\l(t)^n\|h\|_{\mathcal{L}}\ \forall\ |t|<\mathcal{E},\ n\ge 1,\ h\in \mathcal{L};
\end{align*}
where $\forall |t|<\mathcal{E}$, $N(t)$ is a projection onto a one-dimensional subspace.
\f{\rm (iii)} If $\s_\psi:=\lim_{n\to\infty}\tfrac1n\Bbb E(\psi_n^2)>0$ as in \eqref{Leonov}, then $\l(t)=1-\s_\psi^2t^2+o(t^2)$ as $t\to 0$ and the {\tt conditional central limit theorem} holds:
\begin{align*}\tag{{\tt con-CLT}}\label{con-CLT}\T^n1_{[\frac{\psi_n}{\s\sqrt n}\le t]}\xrightarrow[n\to\infty]{}\tfrac1{\sqrt{2\pi}}\int_{-\infty}^te^{-\frac{s^2}2}ds;
\end{align*}
\endproclaim
\f{\bf Remarks}

\Par1\ If, in addition, $P_t$ is a Doeblin-Fortet operator $\forall\ t\in\Bbb R$ and $\psi$ is $T$-{\it aperiodic} in the sense that 
$e^{it\psi}=\l g\circ \tau/g$ with $t\in\Bbb R,\ \l\in\Bbb C,\ |\l|=1\ \&\ g:\Bbb T\to\Bbb C$ entails $t=0,\ \l=1\ \&\ g$ constant, then the conditional local limit theorem holds:
for $I\subset\Bbb R$  an interval, and $k_n\in\Bbb Z,\ \tfrac{k_n}{\s\sqrt n}\to\kappa\in\Bbb R$:
\begin{align*}\tag{{\tt con-LLT}}\label{con-LLT}
 \s\sqrt n\T^n(1_{[\Psi_n\in n\Bbb E(\Psi)+k_n+I]})\xrightarrow[n\to\infty]{} \tfrac{|I|}{\sqrt{2\pi}}e^{-\frac{\kappa^2}2}
\end{align*}
and the skew product $(X\x\Bbb R,m\x\text{\tt Leb},T_\psi)$ is ergodic where $T_\psi(x,y)=Tx,y+\psi(x))$. See \cite{G,G-H,A-De2}.
\

\Par2\ If $(X,m,T)$ is an AFU map as in \cite{Roland} (e.g. an Adler arc map), then for $\psi\in\text{\tt BV}$, the characteristic function operator $P_{t,\psi}$ is a Doeblin-Fortet operator on $(L^1(m),\text{\tt BV})\ \forall\ t\in\Bbb R$, See \cite[\S5]{ADSZ}.

To obtain \eqref{con-CLT} for the stationary process $(X,m,T,\psi)$ via Nagaev's theorem, we verify \eqref{clubsuit} and then $\s_\psi>0$.
\proclaim{Theorem 4.2\ \ \ Analyticity of  Perturbation}
\

Suppose that $\phi:\mathbb{D}\righttoleftarrow$ is non-M\"obius inner with $\phi(0)=0$.
\

Let $b>1$ and let $\u w$ be a summable weight so that $\sum_{n\ge 1}\tfrac{b^n}{w(n)}<\infty$.
\

\ If  $\psi\in \mathcal{h}_{\u w}$ then $\forall\ z\in\Bbb C$,
$\Pi_z:=\Pi_{z,\psi}\in\hom(\mathcal{k}_b,\mathcal{k}_b)$.
\

Moreover  $\Bbb C\ni z\mapsto \Pi_z\in\hom(\mathcal{k}_b,\mathcal{k}_b)$ is holomorphic with
$$\frac{d^n\Pi_z}{d\,z^n}(f)=\Pi_z(\psi^n f)=:\Pi_z^{(n)}(f).$$
In particular \eqref{clubsuit} holds.\endproclaim

\

\

 \proclaim{ Lemma 4.3 (Multiplier lemma)}
\

\f {\rm (i)} Let $b>1$ and let $\u w$ be a summable weight so that $\sum_{n\ge 1}\tfrac{b^n}{w(n)}<\infty$, then
$\mathcal{h}_{\u w}\subseteq M(\mathcal{k}_b)$ with
\begin{align*}\tag*{\dscommercial}\label{dscommercial}
\|f\|_{M(\mathcal{k}_b)}\le R_{b,\u w}\|f\|_{\u w}\ \forall\ f\in \mathcal{h}_{\u w}.
\end{align*}
where $R_{b,\u w}:=\sqrt{1+2\sum_{n\ge 1}\tfrac{b^n}{w(n)}}$.
\

\f{\rm (ii)}\ \  If $\u w$ is a weight so that $\varlimsup_{n\to\infty}\tfrac{w(n+1)}{w(n)}=\infty$, then
        there are no non-constant multipliers of $\mathcal{h}_{\u w}$.
        \endproclaim
     Note that if $w(n)=B^n$ with $B>b$ then $R_{b,\u w}=\sqrt{\tfrac{B+b}{B-b}}$.

        \demo{Proof of (i)}\ \ We have that

\begin{align*}
 \|fg\|_{\mathcal{k}_b}^2&=\sum_{n\in\Bbb Z}b^{|n|}|\widehat{fg}(n)|^2=
\sum_{n\in\Bbb Z}b^{|n|}|\sum_{k\in\Bbb Z}\widehat{f}(k)\widehat{g}(n-k)|^2\\ &=
\sum_{n\in\Bbb Z}b^{|n|}|\sum_{k\in\Bbb Z}\sqrt{w(|k|)}|\widehat{f}(k)|\tfrac{|\widehat{g}(n-k)|}{\sqrt{w(|k|)}}|^2\le
\sum_{n\in\Bbb Z}b^{|n|}\sum_{k\in\Bbb Z}w(|k|)|\widehat{f}(k)|^2\sum_{\ell\in\Bbb Z}\tfrac{|\widehat{g}(n-\ell)|^2}{w(|\ell|)}
\\ &=\|f\|_{\u w}^2\sum_{n,\ell\in\Bbb Z}b^{|n|}\tfrac{|\widehat{g}(n-\ell)|^2}{w(|\ell|)}=
\|f\|_{\u w}^2\sum_{n,\ell\in\Bbb Z}\tfrac{b^{|\ell|}}{w(|\ell|)}b^{|n|-|\ell|}|\widehat{g}(n-\ell)|^2\\ &\le
\|f\|_{\u w}^2\sum_{n,\ell\in\Bbb Z}\tfrac{b^{|\ell|}}{w(|\ell|)}b^{|n-\ell|}|\widehat{g}(n-\ell)|^2=
\|f\|_{\u w}^2\|g\|_{\mathcal{k}_b}^2\sum_{\ell\in\Bbb Z}\tfrac{b^{|\ell|}}{w(|\ell|)}\\ &=R_{b,\u w}^2\|f\|_{\u w}^2\|g\|_{\mathcal{k}_b}^2
\ \ \ \ \ \text{\CheckedBox\ (\eqref{dscommercial})}
\end{align*}
 \demo{Proof of (ii)}\ \ Suppose otherwise, then
        $\exists\ f\in M(\mathcal{h}_{\u w})\ \&\ \ell\ge 1$ with $\widehat{f}(\ell)\ne 0$. Suppose that $\nu_k\uparrow \infty$
        is so
        that $w(\nu_k+1)\ge kw(\nu_k)$ for $k\ge 1$. It follows that
        $$\infty>\|f\|^2_{M(\mathcal{h}_{\u w})}\ge \tfrac{\|f\chi^{\nu_k}\|^2_{\mathcal{h}_{\u w}}}{w(\nu_k)}\ge
        |\widehat{f}(\ell)|^2\tfrac{w(\ell+\nu_k)}{w(\nu_k)}\ge k|\widehat{f}(\ell)|^2
        \xrightarrow[k\to\infty]{}\infty.\ \XBox\ \CheckedBox\ \text{\rm (ii)}$$

\demo{Proof of  Theorem 4.2}
\

It suffices to show that $\forall\ \om\in\Bbb C,\
\exists\ \e=\e_\om>0$ so that
\begin{align*}\tag*{\faShip}\label{faShip}
 \sum_{n=0}^N\tfrac{(z-\om)^n}{n!}\Pi_\om^{(n)}\xrightarrow[N\to\infty]{\hom(\mathcal{k}_b,\mathcal{k}_b)}\Pi_z\ \forall\ z\in B(\om,\e).
\end{align*}

To this end, fix $\om\in\Bbb C$ and $1<b<b_1<b_2$ so that $\sum_{n\ge 1}\tfrac{b_2^n}{w(n)}<\infty$.
\

Let $f\in \mathcal{k}_b$, then for $k\ge 1$,
\begin{align*}\|\Pi^{(k)}_\om(f)\|_{\mathcal{k}_b} &=\|\ttau(\psi^ke^{\om\psi}f)\|_{\mathcal{k}_b}\le M\rho
\|\psi^ke^{\om\psi}f\|_{\mathcal{k}_b}\ \text{by (\eqref{Wheelchair})}\\ &\le M\rho R_{b,{\u w}_{1,b_1}}\|f\|_{\mathcal{k}_b}\|\psi^ke^{\om\psi}\|_{b_1}\ \text{by \eqref{dscommercial}} \\ &=M\rho \sqrt{\tfrac{b_1+b}{b_1-b}}\|f\|_{\mathcal{k}_b}\|\psi^ke^{\om\psi}\|_{b_1}.
\end{align*}
To continue,  by repeated application of  (\eqref{dscommercial}) 
\begin{align*}\|\psi^ke^{\om\psi}\|_{b_1}&\le \left(\sqrt{\tfrac{b_2+b_1}{b_2-b_1}}\right)^{k+1}\|\psi\|_{b_2}^k\|e^{\om\psi}\|_{b_2}\\ &\le
(\sqrt{\tfrac{b_2+b_1}{b_2-b_1}})^{k+1}\|\psi\|_{b_2}^k\exp[|\om|R_{b_2,\u w}\|\psi\|_{\u w}^2].
\end{align*}
Thus
$$\|\Pi_\om^{(k)}\|_{\hom(\mathcal{k}_b,\mathcal{k}_b)}\ll (\sqrt{\tfrac{b_2+b_1}{b_2-b_1}})^k$$ and (\eqref{faShip}) holds  with $\e_\om=\sqrt{\tfrac{b_2-b_1}{b_2+b_1}}$.\ \CheckedBox

\subsection{ Periodicity}
\

\ Let $\phi:\mathbb{D}\righttoleftarrow$ be non-M\"obius inner with $\phi(0)=0$ and let $\psi:\Bbb T\to\Bbb R$ be measurable.
\

We'll call  $t\in\Bbb R$ a $\tau$-{\it period} of $\psi$ if $\exists\ t\in\Bbb R,\ \l\in\Bbb C,\ |\l|=1\ \&\ g:\Bbb T\to\Bbb C$ measurable, so that
$e^{it\psi}=\l g\circ \tau/g$.

\

\

We denote the collection of $\tau$-periods of $\psi$ by $\mathcal{Q}(\psi)$ and call
$\psi$: $\tau$-{\it aperiodic} if $\mathcal{Q}(\psi)=\{0\}$ and $\tau$-{\it periodic} otherwise.

\

It is standard to show that for $t\in\Bbb R,\ \l\in \Bbb C,\ |\l|=1$ and
$f\in L^1(m)$
\begin{align*}\tag*{\faPaw}\label{faPaw}
 e^{it\psi}f=\l f\circ \tau\ \iff\ \ P_t(f):=\ttau(e^{it\psi}f)=\l f
\end{align*}
and also, if $\exists\ \lim_{n\to\infty}\tfrac1n\Bbb E(\psi_n^2)=:\s^2_\psi\ge 0$,
then $\s_\psi>0$ if $\mathcal{Q}(\psi)$ is discrete.

\proclaim{Theorem 4.5 }\
\ \

 Let $\phi:\mathbb{D}\righttoleftarrow$ be non-M\"obius inner with
 with
Denjoy-Wolff point  in $\mathbb{D}$ and nonempty singularity set.

 \f {\rm (i)}\ \ If  $\psi:\Bbb T\to\Bbb R$ is non-constant, real analytic , then
$\mathcal{Q}(\psi)$ is discrete and the stationary process $(\Bbb T,\pi_{\mathfrak{d}_\phi},\tau,\psi)$ satisfies \eqref{con-CLT}.

\f {\rm (ii)}\ \ If both $\mathcal{Z}(\phi)\ \&\ \mathfrak{s}_\phi$ finite, and  $\psi:\Bbb T\to\Bbb R$ is non-constant,
of bounded variation, then $\psi$ is $\tau$-{aperiodic} and the stationary process $(\Bbb T,\pi_{\mathfrak{d}_\phi},\tau,\psi)$ satisfies \eqref{con-LLT}.

\endproclaim
The assumption $\mathfrak{s}_\phi\ne\emptyset$ is essential.
If $\phi$ is a finite Blaschke product  with $\mathfrak{d}(\phi)\in\mathbb{D}\ \ \&\ g:\Bbb T\to\Bbb R$
non-constant and real analytic, then so is $\psi:=g-g\circ\tau:\Bbb T\to\Bbb R$ whence
$P_t(\chi(tg))=\chi(tg)\ \forall\ t\in\Bbb R$ and $\mathcal{Q}(\psi)=\Bbb R$.
    \demo{Proof of (i)}\
    \

    Let $\mathcal{E}>0$ be as in Nagaev's theorem (i). Fix $1<b<B$ with $\psi\in \mathcal{k}_B$. If $\mathcal{Q}(\psi)$ is not discrete, then
$\exists\ t\in (0,\mathcal{E})\cap\mathcal{Q}(\psi)$ and $f\in L^1(m),\ \l\in\Bbb C,\ |\l|=1$
so that $P_t(f)=\l f$. For $t\in (0,\mathcal{E})\cap\mathcal{Q}(\psi),\ P_t$  is a Doeblin-Fortet operator
on $(L^2(m),\mathcal{k}_b)$.     Write $Q:=\overline{\l} P_t$, then $Q:\mathcal{k}_b\righttoleftarrow$ is a Doeblin-Fortet operator, $Qf=f$ and by Lemma 3.1, $f\in \mathcal{k}_b$.

\

By \eqref{faPaw} we have

\begin{align*}\tag*{\Pointinghand}\label{Pointinghand}
 f(\tau x)=\overline{\l}f(x)e^{it\psi(x)}\ \ \text{for $m$-a.e.}\ x\in\Bbb T,
\end{align*}

Using harmonization, we see that \eqref{Pointinghand} holds $\forall x\in\Lambda_{\phi}$.

Now suppose that $\upsilon$ is a point of discontinuity of $\tau$, then $\chi(\upsilon)\in\frak s_\phi$.

\

By Seidel's theorem (\cite{Seidel}, also \cite[theorem 7.48]{Zygmund}), $\forall\ w\in\Bbb T,\ \exists\ z_n=z_n(w)\in\Lambda_{\phi}$ so that $z_n\to\upsilon$
and so that $\tau(z_n)=w$.
Thus
$$f(w)=f(\tau(z_n))=\overline{\l}f(z_n)\psi(z_n)\xrightarrow[n\to\infty]{}\overline{\l}f(\upsilon)\psi(\upsilon)$$
whence $f$ is constant, $e^{it\psi}\equiv\l$, whence $\psi$ (being continuous) is constant.\ \ \CheckedBox\ (i)
\

\demo{Proof of (ii)}
\

Under the assumptions of (ii), $(\Bbb T,\pi_{\mathfrak{d}(\phi)},\tau)$ is an Adler arc map, whence an AFU map as in \cite{Roland}
and $P_{t,\psi}$ is a Doeblin-Fortet operator on $(L^1(m),\text{\tt\small BV}(\Bbb T))\ \forall\ t\in\Bbb R$.
Using \cite[\S5]{ADSZ}, it suffices to prove that $\psi$ is $\tau$-aperiodic.

To this end, suppose that  $f:\Bbb T\to\Bbb C$ is measurable and satisfies \eqref{Pointinghand}.
By Lemma 3.1, $f\in\text{\tt BV}$. We must show that  $f$ is  constant. Suppose otherwise,  then $\exists\ \xi,\ \zeta\in\Bbb T$ so that $|f(\xi)-f(\zeta)|=:\eta>0$.
\

Let $\upsilon\in \chi^{-1}\mathfrak{s}_\phi$, then $\exists\ z_n\to \upsilon$ monotonically so that $\tau(z_{2n})=\xi\ \&\ \tau(z_{2n+1})=\zeta$.
The function $G\in\text{\tt BV}$ where $G(x):=\overline{\l}f(x)e^{it\psi(x)}$, thus so is $f\circ\tau=G$. Therefore
\begin{align*}\infty> \bigvee f\circ\tau
\ge
\sum_{n\ge 1}|f(\tau(z_n))-f(\tau(z_{n+1}))|
=\sum_{n\ge 1}\eta=\infty.\ \ \XBox\ \ \CheckedBox\ \text{(ii)}
\end{align*}

    \

\section{ A local version of Aleksandrov's theorem}\label{loc-Aleks}
\

Fix $A\in\B(\Bbb T),\ m(A)>0\ \&\ \mathcal{C}\subseteq\B(A)$ be a sub-$\s$-algebra.
\

For $1\le p\le\infty$ write

$$H_0^p(A,\mathcal{C}):=\{f\in H^p_0:\ f|_A \ \text{is $\mathcal{C}$-measurable}\}.$$
By \cite[Theorem 17.18]{Rudin} if $f,\ g\in H_0^2(m)$ and $f|_A\equiv g|_A$, then $f\equiv g$.

Let $P:H_0^2\to H_0^2(A,\mathcal{C})$ be the orthogonal projection in the sense that
\begin{align*}\tag*{\faThumbsOUp}\label{faThumbsOUp}
(\text{\tt Id}-P)H_0^2\subset H_0^2(A,\mathcal{C})^\perp;
\end{align*}

$\&$  call  $(A,\mathcal{C})$ an {\it analytic pair} if

$$\mathbb{E}_{m_A}^\mathcal{C}(f)=(Pf)|_A \ \ \ \text{for}\ f\in H_0^2$$
where $m_A(B):=m(A\cap B)\ \&\ \mathbb{E}_{m_A}^\mathcal{C}(f)$ is conditional expectation on the measure space $(A,m_A)$ with respect to the sub-$\s$-algebra $\mathcal{C}\subseteq\B(A)$.
\

\

\

\par  Let $A\in\B(\Bbb T),\ m(A)>0$. It is easy to see that $(A,\mathcal{C})$ is an analytic pair if
either $\mathcal{C}\overset{m_A}=\B(A)$ (in which case $H_0^2(A,\mathcal{C})=H_0^2$); or $\mathcal{C}\overset{m_A}=\{\emptyset,A\}$ (in which case $H_0^2(A,\mathcal{C})=\{0\}$).
\

Next, we give an example which turns out to be general.
\

\proclaim{\bf Example 5.1}\ \  If $A\in\tau^{-1}\B(\Bbb T)\ \&\ \ \mathcal{C}=\tau^{-1}\B(\Bbb T)\cap A$
where $\tau=\tau(\phi)$ with $\phi:\mathcal{D}\righttoleftarrow,\ \phi(0)=0$ inner,\  then $(A,\mathcal{C})$ is a  analytic pair.\endproclaim\pf
\

We have
\begin{align*}
H_0^2(A,\tau^{-1}\B)&:=\{f\in H^2_0:\ f|_A \ \text{is $\tau^{-1}\B$-measurable}\}\\ &=
\{f\in H_0^2:\ \exists\ g\in H_0^2,\ f|_A=g\circ\tau|_A\}\\ &=H_0^2\circ\tau
\ \ \ \  \text{by\ \cite[Theorem 17.18]{Rudin}}.
\end{align*}
To continue, 
note that  because\ $m\circ\tau^{-1}=m$,
$$E^{\tau^{-1}\B(\Bbb T)}(f)=\ttau(f)\circ\tau.$$
\

Since $\ttau:H_0^2\righttoleftarrow$,
$$P:=E_m^{\tau^{-1}\B(\Bbb T)}=\ttau(f)\circ\tau:H_0^2\to H_0^2\circ\tau$$
is orthogonal projection.

\

Let $A=\tau^{-1}B$,  then
$$\mathbb{E}_{m_A}^\mathcal{C}(f)=\mathbb{E}_m^{\tau^{-1}\B}(1_B\circ\tau f)|_A=(1_B\ttau f)\circ\tau)|_A=(Pf)|_A.\ \ \CheckedBox$$
\proclaim{Theorem 5.2}\label{thm5.2}\ \ If $(A,\mathcal{C})$ is a  analytic pair, then either $\mathcal{C}=\{\emptyset,A\}$, or  $\exists\ \phi:\mathcal{D}\righttoleftarrow$  inner with $\phi(0)=0$  so that $A\in\tau^{-1}\B$ and
$\mathcal{C}=A\cap\tau^{-1}\B$ where $\tau=\tau(\phi)$.\endproclaim
The cases with $C=\Bbb T$ are established in \cite{Aleksandrov86}\ \ \footnote{ See also \cite[Theorem 5.6]{Saksman}}.

\demo{Proof} We claim first that

\Par1 $P(gf)=gP(f)\ \forall\ f\in H^\infty_0,\ g\in H^2_0(A,\mathcal{C})$.
\pf
$$P(gf)|_A=\mathbb{E}_{m_A}^\mathcal{C}(fg)=g|_A\mathbb{E}_{m_A}^\mathcal{C}(f)=(gP(f))|_A,$$
and by \cite[theorem 17.18]{Rudin}, $P(gf)=gP(f)$ a.s.  on $\Bbb T$.\  \CheckedBox\P1
\

\Par2 $\forall\ f\in L^2(A,\mathcal{C},m_A)\ \exists\ g,h\in H^2_0(A,\mathcal{C})\ \&\ \g\in\Bbb C$ so that
$$f=(g+\overline{h}+\g)|_A.$$
\pf $\exists\ G,H\in H^2_0\ \&\ \g\in \Bbb C$ so that
$$f=G+\overline{H}+\g.$$
Next, $g=PG,\ h=PH\in H^2_0(A,\mathcal{C})$ and a.s. on $A$:
$$(g+\overline{h}+\g)=(PG+\overline{PH}+\g)=E_{m_A}^\mathcal{C}(f)=f.\ \  \CheckedBox\text{\P2}$$
\

\Par3\ Let $A\in\B(\Bbb T)$ and let $\mathcal{C},\ \mathcal{C}'\subseteq\B(A)$ so that
both $(A,\mathcal{C})$ and  $(A,\mathcal{C}')$ are analytic pairs, then
$$\mathcal{C}\le\mathcal{C}'\ \ \iff\ H^2_0(A,\mathcal{C})\subset H^2_0(A,\mathcal{C}').$$
\demo{\tt Proof} $\Rightarrow$ follows from the definition and  $\Leftarrow$ follows from \P2.\ \CheckedBox
\

\Par4 \ \  $L:=\overline{\text{\tt span}}\{fg:\ f\in H^2_0(A,\mathcal{C}),\ g\in H^\infty_0(A,\mathcal{C})\}\ \ \subsetneq \ \ H^2_0(A,\mathcal{C})$.
\pf
\

Let $$d:=\min\,\{k\ge 0:\ \chi^{-k}g\in H^2\ \forall\ g\in H_0^2(A,\mathcal{C})\},$$ then $d\in\Bbb N$.
\

By \P2, $\exists$\ non constant $G\in H_0^\infty(A,\mathcal{C})$ with $\chi^{-1}G\in H^2$ with the consequence that
$\chi^{-d-1}g\in H^2\ \forall\ g\in L$ whence  $H^2_0(A,\mathcal{C})\setminus L\ne \emptyset$.\ \CheckedBox
\

\Par5 Any $\phi\in H^2_0(A,\mathcal{C})\cap L^\perp,\ \|\phi\|_2=1$ is inner with $\phi(0)=0$.
\pf
\

For $n\ge 1,\ P(\phi\chi^n)\overset{\text{\tiny\P1}}=\phi P(\chi^n)\in L$ whence $\phi\perp\phi P(\chi^n)$    ($\because\ \phi\in L^\perp$).
\

Thus for $n\ge 1$,
\begin{align*}\widehat{(|\phi|^2)}(n)=\<\phi,\phi\chi^n\>=
\<\phi,P(\phi\chi^n)\>=\<\phi,\phi P(\chi^n)\>=0,
 \end{align*}
 whence also $\widehat{(|\phi|^2)} (-n)=\overline{\widehat{(|\phi|^2)} (n)}=0$ and
 $$|\phi|^2=\widehat{(|\phi|^2)} (0)=\|\phi\|_2^2=1.\ \ \CheckedBox\ \P5$$

 \

 Fix $\phi\in H^2_0(A,\mathcal{C})\cap L^\perp$ (inner with $\phi(0)=0$).
 \

 \Par6 \ \ $H^2_0(A,\mathcal{C})=H_0^2\circ\tau=H^2_0(\Bbb T,\tau^{-1}\B)$   where $\tau=\tau(\phi)$.
 \pf Since $\phi\in H_0^\infty(A,\mathcal{C})$, we have
 \par (i) $A\cap\tau^{-1}\B(\Bbb T)\subseteq
 \mathcal{C}$ and (ii) $\phi^n\in H_0^\infty(A,\mathcal{C})\ \forall\ n\ge 1$.
 \

 It follows that any $F\in H_0^2\circ\tau$ is in $H_0^2(A,\mathcal{C})$ having the form
 $F\circ\chi=\sum_{n\ge 1}a_n\phi^n$ with $(a_k:\ k\ge 1)\in\ell^2$.
 \

 Thus
 $H^2_0(A,\mathcal{C})\supseteq H_0^2\circ\tau$.
 \

 To show equality we'll prove that
 $$M:=H^2_0(A,\mathcal{C})\cap (H_0^2\circ\tau)^\perp=\{0\}.$$
 To this end, we show first that
\begin{align*}\tag*{\faAnchor}\label{faAnchor}
\overline{\phi}M\subset M.
\end{align*}
 \pf Let $g\in M$, then $g\perp\phi^j\ \forall\ j\ge 1$ and
  \begin{align*}\tag*{\dsrailways}\label{dsrailways}
    \<\overline{\phi}g,\phi^j\>=0\ \forall\ j\ge 0.
  \end{align*}

 For $k\ge 1$,
 \begin{align*}\<\overline{\phi}g,\chi^{-k}\>&=\<g\chi^k,\phi\>\overset{\text{\tiny\eqref{faThumbsOUp}}}=\<P(g\chi^k),\phi\>=\<gP(\chi^k),\phi\>=0.
 \end{align*}
Thus $\overline{\phi}g\in H^2_0(A,\mathcal{C})$.
\

Each $H\in H^2_0(A,\tau^{-1}\mathcal{B})$ is of form $H=h\circ\tau$ with $h\in H^2_0$. By \eqref{dsrailways},
$$\overline{\phi}g\perp \sum_{k\ge 0}\hat{h}(k)\phi^k=h\circ\tau=H$$
and $\overline{\phi}g\in M$.\ \CheckedBox\ \eqref{faAnchor}
\

To see that $M=\{0\}$ suppose otherwise: that $g\in M,\ g\nequiv 0$,

then by \eqref{faAnchor} (repeatedly) $\overline{\phi}^jg\in M\ \forall\ j\ge 1$ which is impossible unless $g\equiv 0$.\ \CheckedBox\ \P6

\

\Par7 $A\in\tau^{-1}\B(\Bbb T)$,
\pf    By \P6 $H_0^2(A,\mathcal{C})=H_0^2\circ\tau$ and $Pf=\ttau(f)\circ\tau$ and
$$\Bbb E_{m_A}^\mathcal{C}(f)=P(1_Af)|_A$$
$\forall\ f\in H_0^2$ and hence $\forall\ f\in  L^2$.
\

In particular $(\ttau(1_A)\circ\tau)|_A=1$ and we claim that
$\ttau(1_A)\circ\tau=1_A$.
\

To see this  note that
$$\ttau(1_A)\circ\tau=1_A+J$$
where $J:=\ttau(1_A)\circ\tau\cdot 1_{A^c}\ge 0$.
\

Now,
$$m(A)=\Bbb E_m(\ttau(1_A)\circ\tau)=m(A)+\Bbb E_m(J)$$
whence $\Bbb E_m(J)=0,\ J=0$ a.s.,  $\ttau(1_A)\circ\tau=1_A$
and $A\in\tau^{-1}\B(\Bbb T)$.\ \CheckedBox\ \P7


\begin{thebibliography}{ADSZ04}

\bibitem[Aar78]{A-inner78}
Jon Aaronson.
\newblock Ergodic theory for inner functions of the upper half plane.
\newblock {\em Ann. Inst. H. Poincar\'{e} Sect. B (N.S.)}, 14(3):233--253,
  1978.

\bibitem[Aar97]{A1}
Jon Aaronson.
\newblock {\em An introduction to infinite ergodic theory}, volume~50 of {\em
  Mathematical Surveys and Monographs}.
\newblock American Mathematical Society, Providence, RI, 1997.

\bibitem[AD01]{A-De2}
Jon Aaronson and Manfred Denker.
\newblock Local limit theorems for partial sums of stationary sequences
  generated by {G}ibbs-{M}arkov maps.
\newblock {\em Stoch. Dyn.}, 1(2):193--237, 2001.

\bibitem[Adl73]{Fexp}
Roy~L. Adler.
\newblock {$F$}-expansions revisited.
\newblock In {\em Recent advances in topological dynamics ({P}roc. {C}onf.
  {T}opological {D}ynamics, {Y}ale {U}niv., {N}ew {H}aven, {C}onn., 1972; in
  honor of {G}ustav {A}rnold {H}edlund)}, Lecture Notes in Math., Vol. 318,
  pages 1--5. Springer, Berlin, 1973.

\bibitem[ADSZ04]{ADSZ}
J.~Aaronson, M.~Denker, O.~Sarig, and R.~Zweim\"{u}ller.
\newblock Aperiodicity of cocycles and conditional local limit theorems.
\newblock {\em Stoch. Dyn.}, 4(1):31--62, 2004.

\bibitem[Ale86]{Aleksandrov86}
A.~B. Aleksandrov.
\newblock Measurable partitions of the circle induced by inner functions.
\newblock {\em Zap. Nauchn. Sem. Leningrad. Otdel. Mat. Inst. Steklov. (LOMI)},
  149(Issled. Line\u{\i}n. Teor. Funktsi\u{\i}. XV):103--106, 188, 1986.

\bibitem[Ale87]{Aleksandrov87}
A.~B. Aleksandrov.
\newblock Multiplicity of boundary values of inner functions.
\newblock {\em Izv. Akad. Nauk Armyan. SSR Ser. Mat.}, 22(5):490--503, 515,
  1987.

\bibitem[BCJ23]{BCJ}
Oliver Butterley, Giovanni Canestrari, and Sakshi Jain.
\newblock Discontinuities cause essential spectrum.
\newblock {\em Comm. Math. Phys.}, 398(2):627--653, 2023.

\bibitem[Boo57]{Boo}
George Boole.
\newblock On the comparison of transcendents, with certain applications to the
  theory of definite integrals.
\newblock {\em Philosophical Transactions of the Royal Society of London},
  147:745--803, 1857.

\bibitem[Coh80]{Coh}
Donald~L. Cohn.
\newblock {\em Measure theory}.
\newblock Birkh\"{a}user, Boston, Mass., 1980.

\bibitem[Cra91]{Craizer-E}
M.~Craizer.
\newblock Entropy of inner functions.
\newblock {\em Israel J. Math.}, 74(2-3):129--168, 1991.

\bibitem[Den26]{Denjoy}
A.~Denjoy.
\newblock Sur l'it{\'e}ration des fonctions analytiques.
\newblock {\em C. R. Acad. Sci., Paris}, 182:255--257, 1926.

\bibitem[DM78]{De-Meyer}
Claude Dellacherie and Paul-Andr\'{e} Meyer.
\newblock {\em Probabilities and potential}, volume~29 of {\em North-Holland
  Mathematics Studies}.
\newblock North-Holland Publishing Co., Amsterdam-New York, 1978.

\bibitem[DM91]{Man-Do}
Claus~I. {Doering} and Ricardo {Ma\~n\'e}.
\newblock {\em {The dynamics of inner functions}}, volume~3.
\newblock Rio de Janeiro: Sociedade Brasileira de Matem\'atica, 1991.

\bibitem[Don65]{Donoghue}
William~F. Donoghue, Jr.
\newblock On the perturbation of spectra.
\newblock {\em Comm. Pure Appl. Math.}, 18:559--579, 1965.

\bibitem[GH88]{G-H}
Y.~Guivarc'h and J.~Hardy.
\newblock Th\'eor\`emes limites pour une classe de cha\^\i nes de {M}arkov et
  applications aux diff\'eomorphismes d'{A}nosov.
\newblock {\em Ann. Inst. H. Poincar\'e Probab. Statist.}, 24(1):73--98, 1988.

\bibitem[Gla77]{Glaisher}
J.~W.~L. Glaisher.
\newblock Note on certain theorems in definite integration.
\newblock Messenger (2) {VII}. 63-74 (1877)., 1877.

\bibitem[Gor04]{Gordin-Martingale}
M.~I. Gordin.
\newblock A remark on the martingale approximation method for proving the
  central limit theorem for stationary random sequences.
\newblock {\em Zap. Nauchn. Semin. POMI}, 311:124--132, 299--300, 2004.

\bibitem[Gui89]{G}
Y.~Guivarc'h.
\newblock Propri\'et\'es ergodiques, en mesure infinie, de certains syst\`emes
  dynamiques fibr\'es.
\newblock {\em Ergodic Theory Dynam. Systems}, 9(3):433--453, 1989.

\bibitem[Hei77]{Hei}
Maurice Heins.
\newblock On the finite angular derivatives of an analytic function mapping the
  open unit disk into itself.
\newblock {\em J. London Math. Soc. (2)}, 15(2):239--254, 1977.

\bibitem[HH01]{HH}
Hubert Hennion and Lo{\"{\i}}c Herv{\'e}.
\newblock {\em Limit theorems for {M}arkov chains and stochastic properties of
  dynamical systems by quasi-compactness}, volume 1766 of {\em Lecture Notes in
  Mathematics}.
\newblock Springer-Verlag, Berlin, 2001.

\bibitem[ITM50]{IT-M}
C.~T. Ionescu~Tulcea and G.~Marinescu.
\newblock Th\'{e}orie ergodique pour des classes d'op\'{e}rations non
  compl\`etement continues.
\newblock {\em Ann. of Math. (2)}, 52:140--147, 1950.

\bibitem[IU23]{OlegMar}
Oleg Ivrii and Mariusz Urbański.
\newblock Inner functions, composition operators, symbolic dynamics and
  thermodynamic formalism.
\newblock {\em arXiv}, 2308.16063, 2023.

\bibitem[Jan41]{Jankoff}
W.~Jankoff.
\newblock Sur l'uniformisation des ensembles {{(A)}}.
\newblock {\em C. R. (Dokl.) Acad. Sci. URSS, n. Ser.}, 30:597--598, 1941.

\bibitem[Kat04]{Izzy}
Yitzhak Katznelson.
\newblock {\em An introduction to harmonic analysis}.
\newblock Cambridge Mathematical Library. Cambridge University Press,
  Cambridge, third edition, 2004.

\bibitem[Koe84]{Koenigs}
G.~Koenigs.
\newblock Investigations of the integrals of certain functional equations.
\newblock {\em Ann. Sci. {\'E}c. Norm. Sup{\'e}r. (3)}, 1:1--41, 1884.

\bibitem[Leo61]{Leo}
V.~P. Leonov.
\newblock On the dispersion of time means of a stationary stochastic process.
\newblock {\em Teor. Verojatnost. i Primenen.}, 6:93--101, 1961.

\bibitem[Let77]{Let}
G\'{e}rard Letac.
\newblock Which functions preserve {C}auchy laws?
\newblock {\em Proc. Amer. Math. Soc.}, 67(2):277--286, 1977.

\bibitem[Lus30]{Lusin}
N.~Lusin.
\newblock Sur le probl{\`e}me de {M}. {Jacques} {Hadamard} d'uniformisation des
  ensembles.
\newblock {\em Mathematica, Cluj}, 4:54--66, 1930.

\bibitem[LY73]{L-Y}
A.~Lasota and James~A. Yorke.
\newblock On the existence of invariant measures for piecewise monotonic
  transformations.
\newblock {\em Trans. Amer. Math. Soc.}, 186:481--488 (1974), 1973.

\bibitem[Mar89]{Martin-rest}
N.~F.~G. Martin.
\newblock On ergodic properties of restrictions of inner functions.
\newblock {\em Ergodic Theory Dynam. Systems}, 9(1):137--151, 1989.

\bibitem[Nad81]{Nad-inner}
M.~G. Nadkarni.
\newblock Some measure-theoretic remarks on inner functions and related
  results.
\newblock {\em Sankhy\={a} Ser. A}, 43(3):251--259, 1981.

\bibitem[Nag57]{Nag}
S.~V. Nagaev.
\newblock Some limit theorems for stationary {M}arkov chains.
\newblock {\em Teor. Veroyatnost. i Primenen.}, 2:389--416, 1957.

\bibitem[Neu78]{Neu}
J.~H. Neuwirth.
\newblock Ergodicity of some mappings of the circle and the line.
\newblock {\em Israel J. Math.}, 31(3-4):359--367, 1978.

\bibitem[Nor68]{Nor}
Eric~A. Nordgren.
\newblock Composition operators.
\newblock {\em Canadian J. Math.}, 20:442--449, 1968.

\bibitem[Nor72]{Norman}
M.~Frank Norman.
\newblock {\em Markov processes and learning models}.
\newblock Academic Press, New York-London, 1972.
\newblock Mathematics in Science and Engineering, Vol. 84.

\bibitem[NSiG22]{N-S}
Artur Nicolau and Od\'{\i} Soler~i Gibert.
\newblock A central limit theorem for inner functions.
\newblock {\em Adv. Math.}, 401:Paper No. 108318, 39, 2022.

\bibitem[PP90]{Par-Pol}
William Parry and Mark Pollicott.
\newblock Zeta functions and the periodic orbit structure of hyperbolic
  dynamics.
\newblock {\em Ast\'{e}risque}, (187-188):268, 1990.

\bibitem[RE83]{Rou}
J.~Rousseau-Egele.
\newblock Un th\'eor\`eme de la limite locale pour une classe de
  transformations dilatantes et monotones par morceaux.
\newblock {\em Ann. Probab.}, 11(3):772--788, 1983.

\bibitem[Rok61]{Rokhlin}
V.~A. Rokhlin.
\newblock Exact endomorphisms of a {L}ebesgue space.
\newblock {\em Izv. Akad. Nauk SSSR Ser. Mat.}, 25:499--530, 1961.
\newblock English translation in: Amer. Math. Soc. Transl. (2) 39 (1964), 1-36.

\bibitem[Rud74]{Rudin}
Walter Rudin.
\newblock {\em Real and complex analysis}.
\newblock McGraw-Hill Book Co., New York-D\"{u}sseldorf-Johannesburg, second
  edition, 1974.
\newblock McGraw-Hill Series in Higher Mathematics.

\bibitem[Ryc83]{Rychlik}
Marek Rychlik.
\newblock Bounded variation and invariant measures.
\newblock {\em Studia Math.}, 76(1):69--80, 1983.

\bibitem[Sak07]{Saksman}
Eero Saksman.
\newblock An elementary introduction to {C}lark measures.
\newblock In {\em Topics in complex analysis and operator theory}, pages
  85--136. Univ. M\'{a}laga, M\'{a}laga, 2007.

\bibitem[Sei34]{Seidel}
Wladimir Seidel.
\newblock On the distribution of values of bounded analytic functions.
\newblock {\em Trans. Amer. Math. Soc.}, 36(1):201--226, 1934.

\bibitem[Sha93]{Shapiro}
Joel~H. Shapiro.
\newblock {\em Composition operators and classical function theory}.
\newblock Universitext: Tracts in Mathematics. Springer-Verlag, New York, 1993.

\bibitem[Smi29]{Smirnov}
V.~Smirnoff.
\newblock Sur les valeurs limites des fonctions, regulieres a l'interieur d'un
  cercle.
\newblock Journ. Soc Phys-Math Leningrad 2 (1929), 22-37, 1929.
\newblock English translation in: Topics in interpolation theory, Birkhauser
  (1997), 481-494.

\bibitem[Sri98]{Srivastava}
S.~M. Srivastava.
\newblock {\em A course on {B}orel sets}, volume 180 of {\em Graduate Texts in
  Mathematics}.
\newblock Springer-Verlag, New York, 1998.

\bibitem[vN49]{von-N-section}
John von Neumann.
\newblock On rings of operators. {R}eduction theory.
\newblock {\em Ann. of Math. (2)}, 50:401--485, 1949.

\bibitem[Wol26]{Wolff-I}
J.~Wolff.
\newblock Sur l'it{\'e}ration des fonctions born{\'e}es.
\newblock {\em C. R. Acad. Sci., Paris}, 182:200--201, 1926.

\bibitem[YK41]{Yosida-Kakutani}
K\^{o}saku Yosida and Shizuo Kakutani.
\newblock Operator-theoretical treatment of {M}arkoff's process and mean
  ergodic theorem.
\newblock {\em Ann. of Math. (2)}, 42:188--228, 1941.

\bibitem[Zwe98]{Roland}
Roland Zweim\"uller.
\newblock Ergodic structure and invariant densities of non-{M}arkovian interval
  maps with indifferent fixed points.
\newblock {\em Nonlinearity}, 11(5):1263--1276, 1998.

\bibitem[Zyg02]{Zygmund}
A.~Zygmund.
\newblock {\em Trigonometric series. {V}ol. {I}, {II}}.
\newblock Cambridge Mathematical Library. Cambridge University Press,
  Cambridge, third edition, 2002.
\newblock With a foreword by Robert A. Fefferman.

\end{thebibliography}
\end{document}